\documentclass{amsart}

\usepackage{amsfonts}
\usepackage{amssymb}
\usepackage{latexsym}

\newtheorem{theorem}{Theorem}[section]
\newtheorem{lemma}[theorem]{Lemma}

\newtheorem{definition}[theorem]{Definition}
\newtheorem{example}[theorem]{Example}

\newtheorem{claim}[theorem]{Claim}

\theoremstyle{remark}

\numberwithin{equation}{section}

\newcommand{\C}{\mathbb{C}}

\newcommand{\R}{\mathbb{R}}
\newcommand{\Z}{\mathbb{Z}}
\newcommand{\F}{\mathbb{F}}
\newcommand{\Imt}{\mbox{Im}\,}

\begin{document}

\title
{Nonvanishing vector fields on orbifolds}

\author{Carla Farsi}
\address{Department of Mathematics, University of Colorado at Boulder, Campus
Box 395, Boulder, CO 80309-0395 } \email{farsi@euclid.colorado.edu}

\author{Christopher Seaton}
\thanks{The second author was partially supported by a Rhodes College Faculty
Development Endowment Grant.}
\address{Department of Mathematics and Computer Science,
Rhodes College, 2000 N. Parkway, Memphis, TN 38112}
\email{seatonc@rhodes.edu}



\begin{abstract}

We introduce a complete obstruction to the existence
of nonvanishing vector fields on a closed orbifold $Q$.
Motivated by the inertia orbifold, the space of multi-sectors, and the generalized orbifold Euler characteristics, we construct for each finitely generated group $\Gamma$ an orbifold called the space of $\Gamma$-sectors of $Q$.  The obstruction occurs as the Euler-Satake characteristics of the $\Gamma$-sectors for an appropriate choice of $\Gamma$; in the case that $Q$ is oriented, this obstruction is expressed as a cohomology class, the $\Gamma$-Euler-Satake class.  We also acquire a complete obstruction in the case that $Q$ is compact
with boundary and in the case that $Q$ is an open suborbifold of a
closed orbifold.

\end{abstract}

\maketitle


\section{Introduction}
\label{sec-intro}

If $M$ is a closed manifold, then it is well-known that $M$ admits a
smooth, nonvanishing vector field if and only the Euler
characteristic of $M$ vanishes (see e.g. \cite{gp}).  For the case
of a closed orbifold $Q$, the fact that the existence of a
nonvanishing vector field ensures the vanishing of the Euler-Satake
characteristic (i.e. Satake's {\bf Euler characteristic as a
$V$-manifold}) is a trivial consequence of Satake's
Poincar\'{e}-Hopf Theorem in \cite{satake2}. In \cite[Corollary
3.4]{seaton1}, the second author offered a different
Poincar\'{e}-Hopf theorem, demonstrating that a nonvanishing vector
field also implies that the Euler characteristic of the underlying
topological space of $Q$ vanishes.  However, the converse of both of
these statements is false; it is easy to construct examples of
orbifolds such that both characteristics vanish, yet whose singular
strata force any vector field to vanish (see \cite{seaton2}).

Similarly, if $M$ is an open manifold or manifold with boundary,
then it is well-known that $M$ always admits a nonvanishing vector
field.  We note that no requirements are made of the behavior of the
vector field on the boundary; i.e. it need not be tangent to the
boundary nor pointing in or out of $M$.  The case of orbifolds is
again not as straightforward, however, as closed components of the
singular strata may force a vector field to vanish.

In \cite{seaton2}, the second author introduced a complete
cohomological obstruction to the existence of nonvanishing vector
fields on closed orbifolds with cyclic local groups.  In this case,
the obstruction was an element of the Chen-Ruan orbifold cohomology
(see \cite{chenruanorbcohom} or \cite{ademleidaruan}), additively
the cohomology of the inertia orbifold.  For cyclic orbifolds, the
cohomology of the inertia orbifold is large enough to produce a
complete obstruction.  Here, we generalize the construction of the
inertia orbifold to introduce the {\bf space of $\Gamma$-sectors} of
a general orbifold $Q$.  Roughly speaking, the inertia orbifold is
the set of pairs $(x, (g))$ where $x$ is an object in an orbifold
groupoid $\mathcal{G}$ presenting $Q$ and $(g)$ the conjugacy class
of an element $g$ in the isotropy group of $x$. Hence, $(g)$ can be
thought of as the conjugacy class of a homomorphism from $\Z$ into
the isotropy group of $x$.  In contrast, the space of
$\Gamma$-sectors is constructed by choosing homomorphisms from a
fixed, finitely generated group $\Gamma$ into the isotropy group.
The orbifold structure of the space of $\Gamma$-sectors is given by
a translation groupoid via an action of the orbifold groupoid
$\mathcal{G}$. When $\Gamma$ is chosen appropriately, the set of
Euler-Satake characteristics of these {\bf $\Gamma$-sectors} acts as
a complete obstruction to the existence of nonvanishing vector
fields on $Q$. When $Q$ is oriented, we define in the cohomology of
the space of $\Gamma$-sectors an Euler class $e_\Gamma^{ES}(Q)$
called the {\bf $\Gamma$-Euler-Satake class of $Q$} that contains
this information. We demonstrate the following.

\begin{theorem}
\label{thrm-closedobstruction}

Let $Q$ be a closed orbifold and $\Gamma$ a finitely generated
group that covers the local groups of $Q$.  Then $Q$ admits a
smooth, nonvanishing vector field if and only if $\chi_{ES}\left(\tilde{Q}_{(\phi)}\right) = 0$
for each $\Gamma$-sector $\tilde{Q}_{(\phi)}$.  In the case that $Q$ is oriented, this is equivalent to
$e_\Gamma^{ES}(Q) = 0$.

\end{theorem}

We will also show Theorems \ref{thrm-boundarycase} and
\ref{thrm-opencase} which give similar results in the case that $Q$
is a compact orbifold with boundary and $Q$ is an open suborbifold
of a closed orbifold, respectively.

The construction of the $\Gamma$-sectors is motivated by a
construction of Tamanoi in \cite{tamanoi1} and \cite{tamanoi2} (see
also \cite{bryanfulman} and \cite{ohmoto}) used to define
generalized orbifold Euler characteristic of a global quotient
orbifold; i.e. an orbifold that admits a presentation as $M/G$ where
$M$ is a smooth manifold and $G$ is a finite group acting smoothly.
Using the techniques of \cite{chenruanorbcohom} and
\cite{ademleidaruan}, we produce a similar construction to orbifolds
that do not admit such a presentation.

Late in the preparation of this paper, the authors became aware of a
similar construction in \cite[pages 4--8]{leida}.  Leida notes (on
page 14) that his space of fixed-point sectors can be identified
with the mapping space of faithful homomorphisms from finite groups
into the orbifold groupoid $\mathcal{G}$; we take this approach
using a fixed, not necessarily finite group, and do not require that
the homomorphisms be faithful.  It is clear that a similar
obstruction theorem can be proven using Leida's construction.
However, the construction contained here relates specifically to and
generalizes existing constructions for quotients, including orbifold
Euler characteristics.  In a forthcoming paper, the authors will
explore the relationship between the $\Gamma$-sectors given here and
other constructions, including the inertia orbifold, the space of
multi-sectors in \cite{chenruanorbcohom} and \cite{ademleidaruan},
and orbifold Euler characteristics.  Here, our focus is the
properties of the construction itself and the obstruction to
nonvanishing vector fields.

The outline of this paper is as follows. In Section
\ref{sec-localdef}, we give the construction of the $\Gamma$-sectors
as well as the definition of the class $e_\Gamma^{ES}(Q)$.  In
Section \ref{sec-topologicalproperties}, we determine the
topological properties of the $\Gamma$-sectors that we will require.
In Section \ref{sec-obstruction}, we prove Theorem
\ref{thrm-closedobstruction}.  We also prove as Theorems
\ref{thrm-boundarycase} and \ref{thrm-opencase}, the analogous
results in the cases of $Q$ compact with boundary and $Q$ an open
suborbifold of a closed orbifold, respectively.

Note that many authors require of an orbifold $Q$ that each local group $G_x$ act within a chart with a fix-point set of codimension at least $2$.  We make this requirement as well.  However, we note that the construction of the $\Gamma$-sectors in Subsections \ref{subsec-localdef} and \ref{subsec-connectedcomponents} does not require this hypothesis.


\section{Definitions}
\label{sec-localdef}


\subsection{The $\Gamma$-Sectors of an Orbifold}
\label{subsec-localdef}

We assume throughout that $Q$ is an $n$-dimensional orbifold; we do
not assume that $Q$ is effective nor admits a presentation as
quotient. We use the definition and notation in \cite{ademleidaruan}
(see also \cite{chenruangwt}, \cite{moerdijkorbgroupintro},
\cite{moerdijkpronk}, and \cite{ruansgt} for background
information). Recall that presentation of $Q$ is given by an
orbifold groupoid ${\mathcal G}$ (a proper \'{e}tale Lie groupoid)
and a homeomorphism between the orbit space $|{\mathcal G}|$ of
${\mathcal G}$ and the underlying space of $Q$ (see \cite[pages
19--23]{ademleidaruan}). We take a fixed orbifold groupoid
${\mathcal G}$ and identify the underlying space of $Q$ with
$|{\mathcal G}|$.  Let $\sigma : G_0 \rightarrow |\mathcal{G}| = Q$
denote the quotient map.  As usual, $G_0$ and $G_1$ denote the space
of objects and arrows, respectively, in ${\mathcal G}$, $s$ and $t$
the source and target map, respectively, and $G_x$ denotes the set
of loops at a point $x \in G_0$, the {\bf isotropy group of $x$}.

Let $p \in Q$ correspond to the $\mathcal{G}$-orbit of $x \in G_0$
so that $\sigma(x) = p$. There is a neighborhood $V_x$ of $x$ in
$G_0$ diffeomorphic to $\R^n$ in such a way that $x$ corresponds to
the origin and the $\mathcal{G}|_{V_x} = (s, t)^{-1}(V_x \times
V_x)$-action corresponds to a linear $G_x$-representation (see
\cite[page 19 and Proposition 1.44, page 21]{ademleidaruan},
\cite[page 8]{moerdijkorbgroupintro}, or \cite[page
15]{moerdijkpronk}). For ease of notation, we will identify $V_x$
with a subset of $\R^n$ without explicit reference to a choice of
diffeomorphism. In this context, we use $\pi_x : V_x \rightarrow Q$
to denote the restriction of the quotient map $\sigma$ to $V_x$ and
$U_p$ to denote $\pi_x(V_x) \subseteq Q$. Then we call $\{ V_x, G_x,
\pi_x \}$ a {\bf linear orbifold chart for $Q$ at $x$}. Whenever we
use this notation, we will assume that the chart has these
properties. In particular, such a chart defines a groupoid
homomorphism $\xi_x : \mathcal{G}|_{V_x} \rightarrow G_x$ where the
group $G_x$ is treated as a groupoid with space of objects $\{ x
\}$, identifying $\mathcal{G}|_{V_x}$ with $G_x \ltimes V_x$. As
$G_x$ acts on $V_x$, for each $y \in V_x$, $\xi_x$ defines a
bijection between $s^{-1}(y) \cap t^{-1}(V_x)$ and $G_x$. In
particular, $\xi_x$ restricts to an injective group homomorphism
denoted $\xi_x^y = (\xi_x)|_{G_y} : G_y \rightarrow G_x$.

If $x^\prime \in G_0$ is another point in the orbit of $x$ so that
$\sigma(x^\prime) = p$, then there is a $g \in G_1$ with $s(g) = x$,
$t(g) = x^\prime$.  By shrinking $V_x$ if necessary, we can assume
that $s$ restricts to a diffeomorphism $s_g$ from a neighborhood of
$g$ to $V_x$.  Then $gG_x g^{-1} = G_{x^\prime}$ so that $G_x$ and
$G_{x^\prime}$ are isomorphic, and $t \circ s_g^{-1}(V_x)$ is a
neighborhood of $x^\prime$ in $G_0$ diffeomorphic to $\R^n$ where
$x^\prime$ corresponds to the origin; hence, $g$ induces an
isomorphism of orbifold charts, and $\pi_{x^\prime} \circ t\circ
s_g^{-1} = \pi_x$.  Hence, we may refer to the isotropy group of a point $p \in Q$.  By this, we of course mean the isotropy group of an $x \in G_0$ with $\sigma(x) = p$, which (up to isomorphism) does not depend on the choice of $x$.

The following definition follows \cite[pages 52--4]{ademleidaruan}.

\begin{definition}[Space of objects of the $\Gamma$-sectors]

Let $Q$ be an orbifold and $\Gamma$ a finitely generated group.  We let ${\mathcal S}_{\mathcal
G}^\Gamma$ denote the set
\[
    {\mathcal S}_{\mathcal G}^\Gamma    =       \bigcup\limits_{x
    \in G_0} \mbox{HOM}(\Gamma, G_x).
\]
A point in ${\mathcal S}_{\mathcal G}^\Gamma$ will be denoted by
$\phi_x$ where $\phi_x \in \mbox{HOM}(\Gamma, G_x)$.
We let $\beta_\Gamma : {\mathcal S}_{\mathcal G}^\Gamma \rightarrow
G_0$ denote the map $\beta_\Gamma (\phi_x) = x$.

\end{definition}

For each $\phi_x \in {\mathcal S}_{\mathcal G}^\Gamma$, pick a
linear chart $\{ V_x, G_x, \pi_x \}$ for $Q$ at $x \in G_0$.  Let
$V_x^{\langle \phi_x \rangle} = \bigcap\limits_{\gamma \in \Gamma}
V_x^{\phi_x(\gamma)}$ denote the fixed-point subset of $\Imt
\phi_x$.  Let the map
\[
    \kappa_{\phi_x} : V_x^{\langle \phi_x \rangle}
    \longrightarrow
    {\mathcal S}_{\mathcal G}^\Gamma
\]
be defined as follows.  For each $y \in V_x$, $\xi_x^y :G_y
\rightarrow G_x$ is an injective group homomorphism.  If $y \in
V_x^{\langle \phi_x \rangle}$, we have of $\Imt \phi_x \leq
\xi_x^y(G_y) \leq G_x$, so that we can define
$\phi_y:=(\xi_x^y)^{-1}\circ\phi_x : \Gamma \rightarrow G_y$. Let
$\kappa_{\phi_x}(y) = \phi_y$.

\begin{lemma}
\label{lem-localhomsasubmanifold}

Let $Q$ be an orbifold and $\Gamma$ a finitely generated group.  The $\left\{ V_x^{\langle \phi_x \rangle}, \kappa_{\phi_x} \right\}$
give ${\mathcal S}_{\mathcal G}^\Gamma$ the structure of a smooth manifold (with connected components of
different dimensions) in such a way that $\beta_\Gamma$ is a smooth
surjective map.

\end{lemma}

\begin{proof}

Fix $\phi_x \in \mathcal{S}_\mathcal{G}^\Gamma$.  It is clear that
$\kappa_{\phi_x}$ is injective, as it is inverted on its image by
$\beta_\Gamma$.  We give ${\mathcal S}_{\mathcal G}^\Gamma$ the
topology induced by the $\kappa_{\phi_x}$.  Then the $\left\{
V_x^{\langle \phi_x \rangle}, \kappa_{\phi_x} \right\}$ define
manifold charts at each point $\phi_x \in {\mathcal S}_{\mathcal
G}^\Gamma$.  If a different linear chart is chosen at $x$, then it
clearly yields an equivalent (manifold) chart for ${\mathcal
S}_{\mathcal G}^\Gamma$.  If $\phi_x, \psi_y \in {\mathcal
S}_{\mathcal G}^\Gamma$ such that $\kappa_{\phi_x}(V_x^{\langle
\phi_x \rangle}) \cap \kappa_{\psi_y}(V_x^{\langle \psi_y \rangle})
\neq \emptyset$, then
\[
    \kappa_{\phi_x}^{-1} \circ \kappa_{\phi_y} : V_y^{\langle \phi_y
    \rangle} \longrightarrow V_x^{\langle \phi_x \rangle}
\]
is a restriction of the associated transition map for the smooth
manifold $G_0$ to a submanifold and hence smooth.  Therefore, the
$\left\{ V_x^{\langle \phi_x \rangle}, \kappa_{\phi_x} \right\}$ and their
finite intersections define an atlas of smooth charts for ${\mathcal
S}_{\mathcal G}^\Gamma$.

For each $\phi_x \in {\mathcal S}_{\mathcal G}^\Gamma$, the map
$\beta_\Gamma \circ \kappa_{\phi_x}$ is the identity on
$V_x^{\langle \phi_x \rangle}$.  Hence $\beta_\Gamma$ is smooth.

\end{proof}

We define a ${\mathcal G}$-action on ${\mathcal S}_{\mathcal
G}^\Gamma$ by letting $g \in G_1$ act via pointwise conjugation. In
other words, if $x = s(g)$, for each $\gamma \in \Gamma$, we set
\begin{equation}
\label{eq-actiononloopspace}
    (g \phi_x) (\gamma) = g (\phi_x(\gamma))g^{-1}.
\end{equation}
Note that $\beta_\Gamma$ is the anchor map of this action, and that
$g : \beta_\Gamma^{-1} (s(g)) \rightarrow \beta_\Gamma^{-1}(t(g))$
as each $g (\phi_x(\gamma))g^{-1}$ is in the isotropy group of
$t(g)$. The properties of a groupoid action follow trivially.  We
let ${\mathcal G}^\Gamma$ denote the action groupoid ${\mathcal G}
\ltimes {\mathcal S}_{\mathcal G}^\Gamma$.

As $\mathcal{G}^\Gamma$ is the action groupoid for a smooth orbifold
groupoid acting on a smooth manifold, $\mathcal{G}^\Gamma$ is an
orbifold groupoid.  Moreover, the anchor map extends to a
homomorphism $\beta_\Gamma : \mathcal{G}^\Gamma \rightarrow
\mathcal{G}$ (see \cite[pages 39--40]{ademleidaruan}).

\begin{definition}[Space of $\Gamma$-sectors of $Q$]

We let $\tilde{Q}_\Gamma$ denote the orbit space of $|{\mathcal
G}^\Gamma|$ along with the orbifold structure given by
$\mathcal{G}^\Gamma$. We call $\tilde{Q}_\Gamma$ the {\bf space of
$\Gamma$-sectors of $Q$}. A point in $\tilde{Q}_\Gamma$ is the
${\mathcal G}$-orbit of a point $\phi_x \in {\mathcal S}_{\mathcal
G}^\Gamma$, denoted ${\mathcal G}\phi_x$.

\end{definition}

Fix $x \in G_0$ with orbit $\sigma(x) = p \in Q$ and pick a linear
chart $\{ V_x, G_x, \pi_x \}$ for $Q$ at $x$.  By Lemma
\ref{lem-localhomsasubmanifold}, for each $\phi_x \in {\mathcal
S}_{\mathcal G}^\Gamma$, $\kappa_{\phi_x} : V_x^{\langle \phi_x
\rangle} \rightarrow {\mathcal S}_{\mathcal G}^\Gamma$ gives a
manifold chart for ${\mathcal S}_{\mathcal G}^\Gamma$ near $\phi_x$.
We let $C_{G_x}(\phi_x)$ denote the centralizer of $\Imt \phi_x$ in
$G_x$ and $\pi_x^{\phi_x}: V_x^{\langle \phi_x \rangle} \rightarrow
V_x^{\langle \phi_x \rangle}/C_{G_x}(\phi_x)$ the quotient map.

Given $g \in s^{-1}(x)$, it is clear from Equation
\ref{eq-actiononloopspace} that $g\phi_x = \phi_x$ if and only if $g
\in C_{G_x}(\phi_x)$.  Hence, the isotropy group of $\phi_x$ in the
groupoid $\mathcal{G}^\Gamma$ is given by $C_{G_x}(\phi_x)$.  Via
$\kappa_{\phi_x}$, the $\mathcal{G}$-action on $\kappa_{\phi_x}
\left(V_x^{\langle \phi_x \rangle} \right) \subseteq
\mathcal{S}_\mathcal{G}^\Gamma$ corresponds to the
$C_{G_x}(\phi_x)$-action on $V_x^{\langle \phi_x \rangle}$. With
this, we have the following.

\begin{lemma}
\label{lem-orbifoldchartsforgammasectors}

Let $Q$ be an orbifold and $\Gamma$ a finitely generated group. For
each $\phi_x \in {\mathcal S}_{\mathcal G}^\Gamma$, the manifold
chart $\left\{ V_x^{\langle \phi_x \rangle}, \kappa_{\phi_x}
\right\}$ for ${\mathcal S}_{\mathcal G}^\Gamma$ near $\phi_x$
induces a linear orbifold chart $\left\{ V_x^{\langle \phi_x
\rangle}, C_{G_x}(\phi_x), \pi_x^{\phi_x} \right\}$ for
$\tilde{Q}_\Gamma$ at $\phi_x$. The restriction
$\mathcal{G}^\Gamma|_{\kappa_{\phi_x}\left(V_x^{\langle \phi_x
\rangle}\right)}$ of $\mathcal{G}^\Gamma$ to
$\kappa_{\phi_x}\left(V_x^{\langle \phi_x \rangle}\right) \subseteq
\mathcal{S}_\mathcal{G}^\Gamma$ is isomorphic as a groupoid to
$C_{G_x}(\phi_x) \ltimes V_x^{\langle \phi_x \rangle}$.

\end{lemma}

Note that $\kappa_{\phi_x}\left(V_x^{\langle \phi_x \rangle}\right)
\subseteq \mathcal{S}_\mathcal{G}^\Gamma$, so that strictly
speaking, we should say that $\left\{
\kappa_{\phi_x}\left(V_x^{\langle \phi_x \rangle}\right),
C_{G_x}(\phi_x), \pi_x^{\phi_x} \right\}$ is a linear orbifold chart
for $\tilde{Q}_\Gamma$.  In this case, we will make explicit use of
the diffeomorphism $\kappa_{\phi_x}$ to avoid confusing
$V_x^{\langle \phi_x \rangle} \subseteq V_x \subseteq G_0$ and
$\kappa_{\phi_x}\left(V_x^{\langle \phi_x \rangle}\right) \subseteq
\mathcal{S}_\mathcal{G}^\Gamma$.

The following is stated for the case of multi-sectors in \cite[page 53]{ademleidaruan}; see also
\cite[page 17]{moerdijkorbgroupintro} and \cite{crainicmoerdijk1}.

\begin{lemma}
\label{lem-moritainvariance}

Let ${\mathcal G}$ be an orbifold groupoid and $\Gamma$ a finitely
generated group.  A homomorphism of groupoids $\Phi : {\mathcal G}
\rightarrow {\mathcal H}$ induces a homomorphism $\Phi_\ast :
{\mathcal G}^\Gamma \rightarrow {\mathcal H}^\Gamma$. If $\Phi$ is a strong equivalence,
then $\Phi_\ast$ is a strong equivalence.

\end{lemma}

\begin{proof}

We let $\Phi_0$ and $\Phi_1$ denote the maps on objects and arrows, respectively, given by the groupoid homomorphism $\Phi$.
To avoid confusion with $\mathcal{G}$, we use the notation
$(\mathcal{G}^\Gamma)_1$ to denote the
space of arrows of $\mathcal{G}^\Gamma$, and $s_{\mathcal{G}^\Gamma}$ and $t_{\mathcal{G}^\Gamma}$ to denote
the respective source and target maps of $\mathcal{G}^\Gamma$.
We use similar notation for $\mathcal{H}$ and $\mathcal{H}^\Gamma$, where $\mathcal{H}$ has space of objects $H_0$, space of arrows $H_1$, etc..
Note that the spaces of objects of $\mathcal{G}^\Gamma$ and $\mathcal{H}^\Gamma$ are $\mathcal{S}_\mathcal{G}^\Gamma$ and $\mathcal{S}_\mathcal{H}^\Gamma$, respectively.

Every element of $(\mathcal{G}^\Gamma)_1$ consists of an arrow $g
\in G_1$ and an object $\phi_x \in \mathcal{S}_\mathcal{G}^\Gamma$
such that $s(g) = x$.  We let $(g, \phi_x)$ denote the corresponding
arrow in $(\mathcal{G}^\Gamma)_1$, so that
$s_{\mathcal{G}^\Gamma}[(g, \phi_x)] = \phi_x$ and
$t_{\mathcal{G}^\Gamma}[(g, \phi_x)] = g\phi_x$.  The map
$(\beta_\Gamma)_1 : (\mathcal{G}^\Gamma)_1 \rightarrow G_1$ is given
by $(g, \phi_x) \mapsto g$.

For each $x \in G_0$, $\Phi_1 : G_1 \rightarrow H_1$
restricts to a group homomorphism from $G_x$ to $H_{\Phi_1(x)}$. We
define a groupoid homomorphism $\Phi_\ast :\mathcal{G}^\Gamma \rightarrow
\mathcal{H}^\Gamma$ as follows.  First, we define the map on objects,
\begin{equation}
\label{eq-definducedhomobjects}
\begin{array}{rccl}
    \Phi_{\ast 0} :& \mathcal{S}_\mathcal{G}^\Gamma      &\longrightarrow    &   \mathcal{S}_\mathcal{H}^\Gamma      \\\\
                :& \phi_x                              &\longmapsto        &   \Phi_1 \circ \phi_x .
\end{array}
\end{equation}
Then $\Phi_{\ast 0}(\phi_x) : \Gamma \rightarrow H_{\Phi_0(x)}$ is a group homomorphism as required.  For each $\gamma \in \Gamma$ and $g \in G_1$
with $s(g) = x$,
\[
\begin{array}{rcl}
    \Phi_{\ast 0}(g\phi_x)(\gamma)
        &=&     \Phi_1 [ g\phi_x(\gamma)g^{-1}]                 \\\\
        &=&     \Phi_1(g) \Phi_1[\phi_x(\gamma)] \Phi_1(g^{-1}) \\\\
        &=&     \Phi_1(g) [\Phi_{\ast 0}(\phi_x)(\gamma)].
\end{array}
\]
Hence,
\[
    \Phi_{\ast 0} (g\phi_x)   =   \Phi_1(g) \Phi_{\ast 0}(\phi_x),
\]
so that $\Phi_{\ast 0}$ is a $\mathcal{G}$-$\mathcal{H}$-equivariant map via $\Phi_1$.

Fixing $\phi_x \in \mathcal{S}_\mathcal{G}^\Gamma$, pick a (manifold) chart $\left\{ V_x^{\langle \phi_x \rangle}, \kappa_{\phi_x} \right\}$
for $\mathcal{S}_\mathcal{G}^\Gamma$ near $\phi_x$ as given by Lemma \ref{lem-localhomsasubmanifold}.  Similarly, pick a linear orbifold chart
$\{ W_{\Phi_0(x)}, (\mathcal{H})_{\Phi_0(x)}, \varpi_{\Phi_0(x)} \}$ for $|\mathcal{H}|$ at $\Phi_0(x)$; by shrinking charts if necessary, we may assume that $\Phi_0(V_x) \subseteq W_{\Phi_0(x)}$.  Then
$\left\{ W_{\Phi_0(x)}^{\langle \Phi_{\ast 0}(\phi_x) \rangle}, \kappa_{\Phi_{\ast 0}(\phi_x)} \right\}$ is a manifold chart for
$\mathcal{S}_\mathcal{H}^\Gamma$ near $\Phi_{\ast 0}(\phi_x)$.  As $\Phi_1$ commutes with each of the structure maps of $\mathcal{G}$ and $\mathcal{H}$,
for each $y \in V_x^{\langle \phi_x \rangle}$,
\[
\begin{array}{rcl}
    \kappa_{\Phi_{\ast 0}(\phi_x)}^{-1} \circ \Phi_{\ast 0} \circ \kappa_{\phi_x}(y)
        &=&     \kappa_{\Phi_{\ast 0}(\phi_x)}^{-1} \circ \Phi_1 \circ (\xi_x^y)^{-1} \circ \phi_x           \\\\
        &=&     \kappa_{\Phi_{\ast 0}(\phi_x)}^{-1} \circ  \left( (\xi_\mathcal{H})_{\Phi_0(x)}^{\Phi_0(y)}\right)^{-1} \circ \Phi_1 \circ \phi_x
                                   \\\\
        &=&     \kappa_{\Phi_{\ast 0}(\phi_x)}^{-1} \circ  \left( (\xi_\mathcal{H})_{\Phi_0(x)}^{\Phi_0(y)}\right)^{-1} \circ \Phi_{\ast 0}(\phi_x)
                                   \\\\
        &=&     \Phi_0(y).
\end{array}
\]
It follows that the map
\[
    V_x^{\langle \phi_x \rangle}    \stackrel{\kappa_{\phi_x}}{\longrightarrow}
    \mathcal{S}_\mathcal{G}^\Gamma  \stackrel{\Phi_{\ast 0}}{\longrightarrow}
    \mathcal{S}_\mathcal{H}^\Gamma  \stackrel{\kappa_{\Phi_{\ast 0}(\phi_x)}^{-1}}{\longrightarrow}
    W_{\Phi_0(x)}^{\langle \Phi_{\ast 0}(\phi_x) \rangle}
\]
is nothing more than the restriction of $\Phi_0$ to $V_x^{\langle \phi_x \rangle} \subseteq G_0$, and hence is smooth.  As this is true for each chart at each $\phi_x \in \mathcal{S}_\mathcal{G}^\Gamma$, $\Phi_{\ast 0}$ is a smooth map.

For each $(g, \phi_x) \in (\mathcal{G}^\Gamma)_1$, we set
\[
    \Phi_{\ast 1}[(g, \phi_x)] = (\Phi_1(g), \Phi_\ast(\phi_x)).
\]
In other words, $\Phi_{\ast 1}[(g, \phi_x)]$ is the arrow in
$(\mathcal{H}^\Gamma)_1$ given by the action of $\Phi_1(g)$ on
$\Phi_\ast(\phi_x)$.  Then\
\[
    s_{\mathcal{H}^\Gamma}\left(\Phi_{\ast 1}[(g, \phi_x)]\right)
    =
    \Phi_\ast(\phi_x)
    =
    \Phi_{\ast 0} \left(s_{\mathcal{G}^\Gamma}[(g, \phi_x)] \right),
\]
and
\[
    t_{\mathcal{H}^\Gamma}\left(\Phi_{\ast 1}[(g, \phi_x)]\right)
    =
    \Phi_\ast(g\phi_x)
    =
    \Phi_{\ast 0} \left(t_{\mathcal{G}^\Gamma}[(g, \phi_x)] \right),
\]
etc., so that $\Phi_{\ast 0}$ and $\Phi_{\ast 1}$ commute with the
structure maps $\mathcal{G}^\Gamma$ and $\mathcal{H}^\Gamma$.  For
each $(g, \phi_x) \in (\mathcal{G}^\Gamma)_1$, as
$s_{\mathcal{G}^\Gamma}$ is a local diffeomorphism, there is a
neighborhood $W_g$ of $g$ in $(\mathcal{G}^\Gamma)_1$ diffeomorphic
to $V_x^{\langle \phi_x \rangle}$.  Via this diffeomorphism and the
corresponding construction for $(\mathcal{H}^\Gamma)_1$, just as in
the case of $\Phi_{\ast 0}$, $\Phi_{\ast 1}$ corresponds to the
restriction of $\Phi_1$ to the submanifold $W_g$ of $G_1$.  It
follows that $\Phi_{\ast 1}$ is smooth, and so $\Phi_\ast$ is a
homomorphism of Lie groupoids.

Now, assume $\Phi$ is a strong equivalence.  Then $\Phi_0 : G_0
\rightarrow H_0$ is a surjective submersion. Moreover, $\Phi_1$
restricts to an isomorphism from $G_x$ to
$(\mathcal{H})_{\Phi_0(x)}$ for each $x \in G_0$. It follows that
for each $x \in G_0$, the map $\Phi_{\ast 0}$ defined in Equation
\ref{eq-definducedhomobjects} is a bijection between
$\mbox{HOM}(\Gamma, G_x)$  and $\mbox{HOM}(\Gamma,
(\mathcal{H})_{\Phi_0(x)})$, so that $\Phi_{0\ast}$ is surjective.
Within local charts, $\Phi_{\ast 0}$ is the restriction of a
surjective submersion; hence, $\Phi_{\ast 0}$ is a surjective
submersion.  Moreover, the identification of $\mathcal{G}$ with the
pullback of $\mathcal{H}$ via $\Phi_0 \times \Phi_0 : G_0 \times G_0
\rightarrow H_0 \times H_0$ identifies the $\mathcal{G}$-action on
$\mathcal{S}_\mathcal{G}^\Gamma$ with the pullback of the
$\mathcal{H}$-action on $\mathcal{S}_\mathcal{H}^\Gamma$ via
$\Phi_{\ast 0} \times \Phi_{\ast 0} : \mathcal{S}_\mathcal{G}^\Gamma
\times \mathcal{S}_\mathcal{G}^\Gamma \rightarrow
\mathcal{S}_\mathcal{H}^\Gamma \times
\mathcal{S}_\mathcal{H}^\Gamma$.  Hence, $\Phi_\ast$ is a strong
equivalence.

\end{proof}

Noting that the Morita equivalence class of an orbifold groupoid is
the same as the Morita equivalence class via strong equivalences
(see \cite[page 21]{ademleidaruan}), it follows that the orbifold
structure of $\tilde{Q}_\Gamma$ depends only on $\Gamma$ and the orbifold
structure of $Q$ and not on the choice of $\mathcal{G}$.


\subsection{Connected Components of the $\Gamma$-Sectors}
\label{subsec-connectedcomponents}

We now parameterize the space of $\Gamma$-sectors following
\cite[page 83]{ademleidaruan}.

\begin{definition}[Equivalence in ${\mathcal G}^\Gamma$]

Let $Q$ be an orbifold and $\Gamma$ a finitely generated group. Let
$\phi_x, \psi_y \in {\mathcal S}_{\mathcal G}^\Gamma$ and suppose
there is a linear chart $\{ V_x, G_x, \pi_x \}$ at $x$ for $Q$ with
$y \in V_x$. We say that {\bf $\psi_x$  locally covers $\phi_y$ with
respect to the linear chart $\{ V_x, G_x, \pi_x \}$}, written
$\phi_x \stackrel{loc}{\curvearrowright} \psi_y$, if there is a $g
\in G_x$ such that $g [(\xi_x^y \circ \psi_y)(\gamma)] g^{-1} =
\phi_x(\gamma)$ for each $\gamma \in \Gamma$.  When we say $\phi_x$
locally covers $\psi_y$, we mean that there exists a linear chart
with respect to which $\phi_x$ locally covers $\psi_y$.

Extending this to an equivalence relation on all of ${\mathcal
G}^\Gamma$, we say two ${\mathcal G}$-orbits of homomorphisms
${\mathcal G}\phi_x$ and ${\mathcal G}\psi_y$ are {\bf equivalent},
written ${\mathcal G}\phi_x \approx {\mathcal G}\psi_y$, if there is
a finite sequence $\phi_{x_0}, \phi_{x_1}, \ldots , \phi_{x_l}$ such
that $\phi_{x_0} \in {\mathcal G}\phi_x$, $\phi_{x_l} \in {\mathcal
G}\psi_y$, and for each $i$, $\phi_{x_i}
\stackrel{loc}{\curvearrowright} \phi_{x_{i+1}}$ or $\phi_{x_{i+1}}
\stackrel{loc}{\curvearrowright} \phi_{x_i}$.  We let $(\phi_x)$
denote the $\approx$-class of ${\mathcal G}\phi_x$; we will refer to
this class simply as $(\phi)$ when there is no specific
representative $\phi_x$ in mind or to emphasize the lack of
dependence on an $x \in G_0$.  We let $T_Q^\Gamma$ denote the set of
$\approx$-classes in $\mathcal{S}_\mathcal{G}^\Gamma$.

\end{definition}

Note that two homomorphisms are equivalent only if they are
connected by a sequence of local coverings in \emph{linear orbifold
charts}. Allowing charts of the form $M/G$ where $M$ is a manifold
and $G$ a finite group results in a different definition.

The following two lemmas allow us to simplify the definition of
$\approx$ when dealing with $\mathcal{G}$-orbits of elements of
$\mathcal{S}_\mathcal{G}^\Gamma$ rather than points in
$\mathcal{S}_\mathcal{G}^\Gamma$ themselves. Lemma
\ref{lem-localconglocalequal} shows that by picking an appropriate
representative of an orbit, the conjugation in the definition of
$\stackrel{loc}{\curvearrowright}$ is unnecessary. Lemma
\ref{lem-localcoveroverorbits} shows that the definition of
$\stackrel{loc}{\curvearrowright}$ is well-defined on orbits as it
can be; i.e. it holds for all elements of an orbit that have
representatives in the same linear chart.

\begin{lemma}
\label{lem-localconglocalequal}

Let $Q$ be an orbifold, $\Gamma$ a finitely generated group, and $\phi_x, \psi_y \in \mathcal{S}_\mathcal{G}^\Gamma$.
If $\phi_x \stackrel{loc}{\curvearrowright} \psi_y$ with respect to
the chart $\{ V_x, G_x, \pi_x \}$, then there is an element
$\psi_{y^\prime}$ of the ${\mathcal G}$-orbit of $\psi_y$ such that
$y^\prime \in V_x$ and
\[
    \xi_x^{y^\prime} \circ \psi_{y^\prime} = \phi_x.
\]

\end{lemma}

\begin{proof}

Suppose there is a linear chart and a $g \in G_x$ such that $g [(\xi_x^y \circ \psi_y)(\gamma)]
g^{-1} = \phi_x(\gamma)$ for each $\gamma \in \Gamma$.
As $\xi_x$ defines a surjective map from $s^{-1}(y) \cap t^{-1}(V_x)$ onto $G_x$, there is an arrow $h \in G_1$
such that $s(h) = y$, $t(h) \in V_x$, and $\xi_x(h) = g$. Therefore,
for each $\gamma \in \Gamma$, recalling that $\xi_x^y$ is simply the restriction of $\xi_x$ to $G_y$,
\[
\begin{array}{rcl}
    \phi_x(\gamma)
        &=&    g [(\xi_x^y \circ \psi_y)(\gamma)] g^{-1}
                                        \\\\
        &=&     \xi_x(h) [(\xi_x^y \circ \psi_y)(\gamma)] \xi_x(h^{-1})
                                        \\\\
        &=&     \xi_x [h (\psi_y(\gamma)) h^{-1}]
                                        \\\\
        &=&     \xi_x^{t(h)}[h (\psi_y(\gamma))
        h^{-1}]                         \\\\
        &&\mbox{(as $h (\psi_y(\gamma)) h^{-1} \in G_{t(h)}$)}
                                        \\\\
        &=&     [\xi_x^{t(h)}(h \psi_y )](\gamma).
\end{array}
\]
Setting $y^\prime = t(h)$ and $\psi_{y^\prime} = h\psi_y$, we are
done.

\end{proof}

\begin{lemma}
\label{lem-localcoveroverorbits}

Let $Q$ be an orbifold, $\Gamma$ a finitely generated group, and $\phi_x, \psi_y \in \mathcal{S}_\mathcal{G}^\Gamma$.
If $\phi_x \stackrel{loc}{\curvearrowright} \psi_y$, then $\phi_x$
locally covers every element of the ${\mathcal G}$-orbit of $\psi_y$
in $\beta_\Gamma^{-1}(V_x)$.

\end{lemma}

\begin{proof}

Suppose $\phi_x \stackrel{loc}{\curvearrowright} \psi_y$ so that there is a $g \in G_x$ such that for each $\gamma \in \Gamma$,
\[
    \phi_x(\gamma)
    =    g [(\xi_x^y \circ \psi_y)(\gamma)] g^{-1}.
\]
Each element of the ${\mathcal G}$-orbit of $\psi_y$ in
$\beta_\Gamma^{-1}(V_x)$ is of the form $h\psi_y$ for some $h \in
G_1$ with $s(h) = y$ and $t(h) \in V_x$.  Fixing one such $h\psi_y$,
we have for each $\gamma \in \Gamma$
\[
\begin{array}{rcl}
    g\xi_x(h^{-1}) \left[\left(\xi_x^{t(h)} \circ (h\psi_y)\right)(\gamma)\right] [g\xi_x(h^{-1})]^{-1}
        &=&     g\xi_x(h^{-1}) \left[\xi_x^{t(h)}  \left(h\psi_y(\gamma)h^{-1}\right)\right] \xi_x(h)g^{-1}
                            \\\\
        &=&     g\xi_x(h^{-1}h) \left[\left(\xi_x^y \circ \psi_y\right)(\gamma)\right]  \xi_x(h^{-1}h)g^{-1}
                            \\\\
        &=&     g [(\xi_x^y \circ \psi_y)(\gamma)] g^{-1}
                            \\\\
        &=&     \phi_x (\gamma).
\end{array}
\]
As $g \xi_x(h^{-1}) \in G_x$, it follows that $\phi_x
\stackrel{loc}{\curvearrowright} h\psi_y$.

\end{proof}

\begin{lemma}
\label{lem-finitegammasectors}

Suppose $Q$ is a compact orbifold and $\Gamma$ is a finitely generated group.  Then
$T_Q^\Gamma$ is finite.

\end{lemma}

\begin{proof}

For each $p \in Q$, pick a linear chart $\{ V_x, G_x, \pi_x \}$ at
$p$.  Then the collection $\{ U_p : p \in Q \}$ forms an open cover
of $Q$.  As $Q$ is compact, pick a finite subcover corresponding to
the points $x_1, x_2, \ldots x_k \in G_0$ with respective orbits
$\sigma(x_i) = p_i \in Q$.  We claim that each $(\psi) \in
T_Q^\Gamma$ has a representative in the set
\[
    {\mathcal H} = \bigcup\limits_{i=1}^k \mbox{HOM}(\Gamma, G_{x_i}),
\]
which is clearly is finite.

Let $\psi_y : \Gamma \rightarrow G_y$ be an arbitrary element of
${\mathcal S}_{\mathcal G}^\Gamma$ and let $q = \sigma(y) \in Q$
denote the orbit of $y$. Then there is an $i$ such that $q \in
U_{p_i}$; hence, there is an $h \in G_1$ with $s(h) = y$ and $t(h)
\in V_{x_i}$. Recalling that $\xi_x^{t(h)} : G_{t(h)} \rightarrow
G_x$ is an injective group homomorphism, define $\phi_x =
\xi_x^{t(h)} \circ (h\psi_y)$, i.e.
\[
    \phi_x : \Gamma  \stackrel{h\psi_y}{\longrightarrow} G_{t(h)}
    \stackrel{\xi_x^{t(h)}}{\longrightarrow} G_x,
\]
and then $\phi_x$ is an element of ${\mathcal H}$.

That $\phi_x \stackrel{loc}{\curvearrowright} h\psi_y$ is obvious
from the definition of $\phi_x$.  Hence, as $h\psi_y \in {\mathcal
G}\psi_y$, this implies that $\phi_x \approx \psi_y$.

\end{proof}

If $\phi_x \stackrel{loc}{\curvearrowright} \phi_y$ with respect to
the linear chart $\{ V_x, G_x, \pi_x \}$, then by Lemma
\ref{lem-localconglocalequal}, the $G_x$-orbit of $y$ in $V_x$
intersects $V_x^{\langle \phi_x \rangle}$.  As $G_x$ acts linearly
so that $V_x^{\langle \phi_x \rangle}$ is a subspace, it follows
that $\phi_x$ and $\phi_y$ represent ${\mathcal G}$-orbits in the
same connected component of $\tilde{Q}_\Gamma$. If $\phi_x \approx
\phi_y$, then they are connected by a finite sequence of points
related by local coverings. For each $i$, the two points
$C_{G_{x_i}}(\phi_{x_i})x_i$ and
$C_{G_{x_{i+1}}}(\phi_{x_{i+1}})x_{i+1}$ lie in the same connected
component of $\tilde{Q}_\Gamma$ regardless of the direction of the
covering, so that the ${\mathcal G}$-orbits of $\phi_x$ and $\phi_y$
lie in the same connected component of $\tilde{Q}_\Gamma$.

Conversely, each chart for $\tilde{Q}_\Gamma$ can be taken to be of
the form $\left\{ V_x^{\langle\phi_x\rangle}, C_{G_x}(\phi_x),
\pi_x^{\phi_x} \right\}$, induced from a chart $\{ V_x, G_x, \pi_x
\}$ at $x$, and the image of the injective homomorphism $\xi_x^y :
G_y \rightarrow G_x$ contains $\Imt \phi_x$ if and only if $y \in
V_x^{\langle \phi_x\rangle}$.  Hence, this chart defines a local
covering by $\phi_x$ of each homomorphism $\phi_y = (\xi_x^y)^{-1}
\circ\phi_x$ corresponding to a $y \in V_x^{\langle \phi_x
\rangle}$.  If $\phi_x$ and $\phi_y$ represent points whose orbits
are in the same connected component of $\tilde{Q}_\Gamma$, then
there is a path connecting the orbits of $\phi_x$ and $\phi_y$. Pick
a linear chart for each point on the path and then a finite subcover
of these uniformized sets.  It follows that there is a finite
sequence of local equivalences connecting $\phi_x$ to $\phi_y$. With
this, we make the following definition.

\begin{definition}[$\Gamma$-sector]

Let $\tilde{Q}_{(\phi)}$ denote the subset of $\tilde{Q}_\Gamma$
corresponding to orbits of points in the $\approx$-class $(\phi)$.  Then
\[
    \tilde{Q}_\Gamma = \bigsqcup\limits_{(\phi) \in T_Q^\Gamma}
    \tilde{Q}_{(\phi)}
\]
is a decomposition of $\tilde{Q}_\Gamma$ into connected components.
We call $\tilde{Q}_{(\phi)}$ the {\bf $\Gamma$-sector corresponding
to $(\phi)$.}

\end{definition}

We let $\pi : \tilde{Q}_\Gamma
\rightarrow Q$ denote the map $\pi(\mathcal{G}\phi_x) = \sigma(x)$
that sends the orbit $\mathcal{G}\phi_x$ in $\tilde{Q}_\Gamma$ to the
orbit of $x$ in $Q$.  Note that $\pi$ is the map on orbit spaces induced by $\beta_\Gamma$ and hence is a smooth map of orbifolds.


\subsection{Euler-Satake Characteristics and the $\Gamma$-Euler-Satake Class}
\label{subsec-eulersatakeclass}

If $Q$ is closed, we will use $\chi_{top}(Q)$ to denote the usual
Euler characteristic of the underlying topological space
$\mathbb{X}_Q$ of $Q$.  We let $\chi_{ES}(Q)$ denote Euler-Satake
characteristic of $Q$. Recall that the Euler-Satake characteristic
of an orbifold $Q$ is a rational number; see \cite{satake2} where
this number is called the {\bf Euler characteristic of $Q$ as a
$V$-manifold}, or \cite{seaton1} where this quantity is
denoted\footnote{We avoid using this notation here, as it is more
frequently used in the literature to denote the stringy orbifold
Euler characteristic of $Q$.} $\chi_{orb}(Q)$.

Recall (see \cite[Definition 2.25 and Definition 2.28, pages 44--45]{ademleidaruan}; see also \cite{chenruangwt}) that an orbifold vector bundle over $Q$ of rank $k$ is given by a $\mathcal{G}$-vector
bundle $\rho : E \rightarrow G_0$ of rank $k$ such that each $g \in G_1$ induces
a linear isomorphism of fibers $g:\rho^{-1}(s(g))\rightarrow
\rho^{-1}(t(g))$.  Sections of this orbibundle correspond to
$\mathcal{G}$-invariant sections $\omega : G_0 \rightarrow E$.  We
denote the translation groupoid $\mathcal{E} = \mathcal{G} \ltimes
E$.  A bundle $E$ is called {\bf good} if for each $x \in G_0$,
$\mbox{Ker}(G_x)$ acts trivially on each fiber $E_x$ where $\mbox{Ker}(G_x)$ denotes the set of constant arrows in $G_x$; i.e. the arrows $g \in G_x$ such that there is a neighborhood of $g$ in $G_1$ on which $s = t$.  Within a linear chart over which $E$ is trivial, this means that the kernel of the group action on the base space coincides with the kernel of the group action on the total space.

Note that as $\mathcal{E}$ is an orbifold structure for
$|\mathcal{E}|$, we can apply the construction of $\Gamma$-sectors
to form $\widetilde{|\mathcal{E}|}_\Gamma = |\mathcal{E}^\Gamma |$
where $\mathcal{E}^\Gamma = \mathcal{E} \ltimes
\mathcal{S}_\mathcal{E}^\Gamma$.

\begin{lemma}
\label{lem-bundlesovergammasectors}

Let $\rho : E \rightarrow G_0$ be an orbifold vector bundle of rank $k$, and let
$\Gamma$ be a finitely generated group.  Then
$\mathcal{S}_\mathcal{E}^\Gamma$ is naturally a
$\mathcal{G}^\Gamma$-vector bundle over
$\mathcal{S}_\mathcal{G}^\Gamma$ making
$\widetilde{|\mathcal{E}|}_\Gamma$ into an orbifold vector bundle over
$\tilde{Q}_\Gamma$. An orientation of $E$ induces an orientation of $\mathcal{S}_\mathcal{E}^\Gamma$.
A $\mathcal{G}$-invariant section $\omega : G_0
\rightarrow E$ naturally induces an $\mathcal{E}^\Gamma$-invariant
section $\tilde{\omega}_\Gamma : \mathcal{S}_\mathcal{G}^\Gamma
\rightarrow \mathcal{S}_\mathcal{E}^\Gamma$.  Within corresponding
linear charts for $\mathcal{E}^\Gamma$ and $\mathcal{E}$,
$\tilde{\omega}_\Gamma$ is simply the restriction of $\omega$ to a
subspace; in particular $\tilde{\omega}_\Gamma (\phi_x) = 0$ if and
only if $\omega(x) = 0$.

\end{lemma}

\begin{proof}

As in proof of Lemma \ref{lem-moritainvariance}, we let $(\mathcal{E})_1$ denote the
space of arrows $\mathcal{E}$, $(\mathcal{E})_e$ the isotropy group of $e \in E$,
and $s_\mathcal{E}$ and $t_\mathcal{E}$ the source and target maps,
respectively.  Note that $E$ is the space of objects of $\mathcal{E}$.
An element of $(\mathcal{E})_1$ is given by a $g \in
G_1$ and an $e \in \rho^{-1}(s(g)) \subseteq E$; we will denote this
arrow $(g, e)$.  For each $e \in E$, the map
\begin{equation}
\label{eq-defnu}
\begin{array}{rccl}
    \rho_1     :&     (\mathcal{E})_1 &\longrightarrow&   G_1
                        \\\\
            :&     (g, e)             &\longmapsto&   g
\end{array}
\end{equation}
restricts to an injective homomorphism from the isotropy group
$(\mathcal{E})_e$ into the isotropy group $G_{\rho(e)}$ of $\rho(e)
\in G_0$.  Regarding the space of $\Gamma$-sectors, the space of
objects in $\mathcal{E}^\Gamma$ is $\mathcal{S}_\mathcal{E}^\Gamma$.
We show that $\mathcal{S}_\mathcal{E}^\Gamma$ is a
$\mathcal{G}^\Gamma$-vector bundle over
$\mathcal{S}_\mathcal{G}^\Gamma$ with the desired properties.

Each element of $\mathcal{S}_\mathcal{E}^\Gamma$ is a homomorphism
$\phi_e : \Gamma \rightarrow (\mathcal{E})_e$ where $e \in E$.  We
define
\[
\begin{array}{rccl}
    \tilde{\rho}_\Gamma
        :& \mathcal{S}_\mathcal{E}^\Gamma  &\longrightarrow&   \mathcal{S}_\mathcal{G}^\Gamma
                        \\\\
        :& \phi_e      &\longmapsto&       \rho_1 \circ \phi_e
\end{array}
\]
Then $\tilde{\rho}_\Gamma(\phi_e): \Gamma \rightarrow G_{\rho(e)}$
is a homomorphism for each $\phi_e \in
\mathcal{S}_\mathcal{E}^\Gamma$.

A linear chart for $|\mathcal{E}^\Gamma|$ whose image contains $\phi_e$ is given as
follows. Let $x = \rho(e) \in G_0$.  Then there is a linear chart
$\{ V_x, G_x, \pi_x \}$ for $Q$ at $x$.  By shrinking $V_x$ if
necessary, we may assume that $E|_{V_x}$ is trivial as a vector
bundle over $V_x$ (although the $G_x$-structure need not be
trivial).  Then $\{V_x \times \R^k, G_x, \tilde{\pi}_x \}$ is a
linear chart for $|\mathcal{E}|$ whose image contains $e$ with $\tilde{\pi}_x(V_x
\times \R^k) = \rho^{-1}(V_x) \subseteq E$. If $\mbox{proj}_{V_x} :
V_x \times \R^k \rightarrow V_x$ denotes the projection onto the
first factor, then $\pi_x \circ \mbox{proj}_{V_x} = \rho \circ
\tilde{\pi}_x$.

Note that we are not interested in linear charts for $\mathcal{E}$
at $e$, as such charts do not respect the structure of $E$ as a rank
$k$ vector bundle over $G_0$.  Specifically, if $e$ does not
correspond to an element of the zero section in $E$, then a chart in
which $e$ corresponds to the identity in $\R^n \times \R^k$ would
not give a local trivialization of $E$.  Rather, $\{V_x \times \R^k,
G_x, \tilde{\pi}_x \}$ is a linear chart whose image contains
complete fibers of $E$ such that $V_x \times\{ 0 \}$ corresponds to
the zero-section, and the origin in $V_x \times\R^n$ corresponds to
$x$ through the injection of $G_0$ into $E$ as the zero-section.

By Lemma \ref{lem-orbifoldchartsforgammasectors}, a chart for
$\widetilde{|\mathcal{E}|}_\Gamma = |\mathcal{E}^\Gamma|$ whose
image contains $\phi_e$ is of the form $\left\{(V_x \times
\R^k)^{\langle\phi_e\rangle}, C_{G_x}(\phi_e),
\tilde{\pi}_x^{\phi_e} \right\}$, and a linear chart for
$\tilde{Q}_\Gamma = |\mathcal{G}^\Gamma|$ at
$\tilde{\rho}_\Gamma(\phi_e)$ is of the form $\left\{V_x^{\langle
\tilde{\rho}_\Gamma(\phi_e) \rangle},
C_{G_x}(\tilde{\rho}_\Gamma(\phi_e)),
\pi_x^{\tilde{\rho}_\Gamma(\phi_e)} \right\}$.  Since $G_x$ acts
linearly on each fiber of $V_x \times \R^k$ and $G_x(V_x \times \{ 0
\}) = (V_x \times \{ 0 \})$, we have
\[
\begin{array}{rcl}
    \left(V_x \times \R^k\right)^{\langle \phi_e\rangle}
        &=&   (V_x \times \{ 0 \})^{\langle \phi_e\rangle} \times (\R^k)^{\langle \phi_e\rangle}
                    \\\\
        &=&     V_x^{\langle \tilde{\rho}_\Gamma(\phi_e) \rangle} \times (\R^k)^{\langle
                \phi_e\rangle}.
\end{array}
\]
By construction,
$\tilde{\pi}_x^{\tilde{\rho}_\Gamma(\phi_e)}\circ\mbox{proj}_{V_x^{\langle
\phi_e \rangle}} = \tilde{\rho}_\Gamma \circ
\tilde{\pi}_x^{\phi_e}$. As such charts exist at each point $\phi_e
\in \mathcal{S}_\mathcal{E}^\Gamma$, and as they clearly transform
appropriately, being restrictions of charts for $E$,
$\mathcal{S}_\mathcal{E}^\Gamma$ is a $\mathcal{G}^\Gamma$-vector
bundle over $\mathcal{S}_\mathcal{G}^\Gamma$.  Moreover, an orientation of $E$ obviously induces an orientation of $\mathcal{S}_\mathcal{E}^\Gamma$ by restriction within charts.

Now, let $(g, \phi_x) \in (\mathcal{G}^\Gamma)_1$ where $g \in G_0$ with $s(g) = x$. Note that
$s_{\mathcal{G}^\Gamma}[(g, \phi_x)] = \phi_x$ is a homomorphism
$\Gamma \rightarrow G_x$. Hence, $(g, \phi_x)$ induces a map
\[
    (g, \phi_x)
    : \tilde{\rho}_{\Gamma}^{-1}(\phi_x)  \longrightarrow
    \tilde{\rho}_{\Gamma}^{-1}(t(\phi_x)).
\]
With respect to a chart for $\mathcal{E}^\Gamma$ as above,
$\tilde{\rho}_{\Gamma}^{-1}(\phi_x) = (\rho^{-1}(x))^{\langle \phi_x
\rangle}$ and $\tilde{\rho}_{\Gamma}^{-1}(t(\phi_x)) =
(\rho^{-1}[t(g)])^{\langle g\phi_x \rangle}$, so that $(g, \phi_x)$
is simply the restriction of the linear isomorphism $g :
\rho^{-1}(s(g)) \rightarrow \rho^{-1}(t(g))$ onto the invariant
subspace $(\rho^{-1}(x))^{\langle \phi_x \rangle}$ with image
$(\rho^{-1}[t(g)])^{\langle g\phi_x \rangle}$. Therefore, it is a
linear isomorphism, and $\tilde{\rho}_\Gamma :
\mathcal{S}_\mathcal{E}^\Gamma \rightarrow
\mathcal{S}_\mathcal{G}^\Gamma$ defines an orbifold vector bundle
$\widetilde{|\mathcal{E}|}_\Gamma \rightarrow \tilde{Q}_\Gamma$.

A section $Q \rightarrow |\mathcal{E}|$ is a $\mathcal{G}$-invariant
section $\omega : G_0 \rightarrow E$.  For each $x \in G_0$, as
$\omega (x)$ is $G_x$-invariant, the map $\rho_1$ defined in
Equation \ref{eq-defnu} maps the isotropy group
$(\mathcal{E})_{\omega(x)}$ of $\omega(x)$ isomorphically onto
$G_x$.  Given a section $\omega$, we define
\[
    \tilde{\omega}_\Gamma
        :     \mathcal{S}_\mathcal{G}^\Gamma
        \longrightarrow   \mathcal{S}_\mathcal{E}^\Gamma
\]
by $[\tilde{\omega}_\Gamma(\phi_x)](\gamma) =
(\rho_1|_{(\mathcal{E})_{\omega(x)}})^{-1}[\phi_x(\gamma)]$ for each
$\gamma \in \Gamma$. In a chart for $\mathcal{E}^\Gamma$ of the form
$\left\{(V_x \times \R^k)^{\langle \tilde{\omega}_\Gamma (\phi_x)
\rangle}, C_{(\mathcal{E})_{\omega(x)}}(\tilde{\omega}_\Gamma
(\phi_x)), \tilde{\pi}_x^{\tilde{\omega}_\Gamma (\phi_x)} \right\}$,
the value of $\tilde{\omega}_\Gamma(y)$ coincides with that of
$\omega(y)$, so that $\tilde{\omega}_\Gamma$ is, within a chart,
just the restriction of $\omega$ to the invariant subspace $(V_x
\times \{ 0 \} )^{\langle \tilde{\omega}_\Gamma (\phi_x) \rangle} =
V_x^{\langle \phi_x\rangle}$.  Hence $\tilde{\omega}_\Gamma$ is a
$\mathcal{G}^\Gamma$-invariant section of
$\mathcal{S}_\mathcal{E}^\Gamma$.  Moreover, as it is locally just a
restriction, $\tilde{\omega}_\Gamma(\phi_x) = 0$ if and only if
$\omega(x) = 0$.

\end{proof}

Note that the bundle $\tilde{\rho}_\Gamma :
\mathcal{S}_\mathcal{E}^\Gamma \rightarrow
\mathcal{S}_\mathcal{G}^\Gamma$ is not simply the pullback of the
bundle $E$ via $\beta_\Gamma$; it generally has different ranks over
different connected components. It is easy to see that the the
operations of forming the tangent bundle, cotangent bundle, and its
exterior powers commute with the operation of forming the
$\Gamma$-sectors; that is,
\[
    T(\tilde{Q}_\Gamma) =   \widetilde{(TQ)}_\Gamma,
\]
\[
    T^\ast(\tilde{Q}_\Gamma) =   \widetilde{(T^\ast Q)}_\Gamma,
\]
etc.  These are, however, all good vector bundles.  We note that it
is possible that the orbifold vector bundle given by $\rho: E
\rightarrow G_0$ is a good vector bundle while the induced bundle on
$\Gamma$-sectors is not. We illustrate this with the following
example.

\begin{example}
\label{ex-goodbundleinducesbadbundle}

Let $G_0 = \C$ and $F = \C^2$ with basis $\{ f_1, f_2 \}$ (note that
we only use the complex structure to simplify notation).  Let $E =
G_0 \times F$ be the trivial bundle and $\rho : E \rightarrow G_0$
the projection.  Let $G = \Z_2 \oplus \Z_3  = \langle \alpha_1
\rangle \oplus \langle \alpha_2\rangle$ act as follows.  On $G_0$,
$\alpha_1$ acts as multiplication by $-1$ and $\alpha_2$ acts as
multiplication by $e^{2\pi i/3}$.  On $E$, we let $\alpha_1(y,
c_1f_1 + c_2f_2) = (\alpha_1y, -c_1 f_1 + c_2 f_2)$ and $\alpha_2(y,
c_1f_1 + c_2f_2) = (\alpha_2y, c_1 f_1 + e^{2\pi i/3} c_2 f_2)$.
Then $E$ defines an orbifold vector bundle over the orbifold $Q =
G\ltimes G_0$.  Note that the kernel of the action on both $G_0$ and
$E$ is trivial, so that $E \rightarrow G_0$ is a good orbifold
vector bundle.

Let $x$ denote the origin in $G_0$ and let $\Gamma = \Z$ be
generated by $\gamma$.  Define $\phi_x : \Z \rightarrow G$ by
$\phi_x : \gamma \mapsto \alpha_1$ and let $(\phi)$ denote the
$\approx$-class of $\phi_x$ as usual.  Then a chart for
$\tilde{Q}_{(\phi)}$ at the orbit of $\phi_x$ is simply the origin
with $\alpha_1$ and $\alpha_2$ acting trivially. The kernel of the
action is all of $G$.

Let $e$ denote the origin in $\rho^{-1}(x)$, and define $\phi_e :
\gamma \mapsto \alpha_1$.  Then $\tilde{\rho}_\Gamma(\phi_e) =
\phi_x$.  A chart for the fiber of
$\widetilde{|\mathcal{E}|}_\Gamma$ at the orbit of $\phi_y$ is
$\mbox{span}_\C(f_2)$ where $\alpha_1$ acts trivially and $\alpha_2$
acts by multiplication by $e^{2\pi i/3}$.  The kernel of the action
is $\langle \alpha_1 \rangle$.

We see that $E \rightarrow G_0$ is a good vector bundle, while the induced bundle over the
space of $\Z$-sectors is a bad vector bundle; on the connected
component $\tilde{Q}_{(\phi)}$, the kernel of the action on the
fiber is a proper subgroup of the kernel of the action on the base
space.

\end{example}

\begin{definition}[$\Gamma$-cohomology]
\label{def-cohomology}

Let $Q$ be an orbifold and $\Gamma$ a finitely generated group. We
let $H_{\Gamma}^\ast(Q)$ denote the (singular or de Rham) cohomology
of $\tilde{Q}_\Gamma$.  If $Q$ admits an almost-complex structure,
we let $H_{orb, \Gamma}^\ast(Q)$ denote Chen-Ruan cohomology of
$\tilde{Q}_\Gamma$ (see \cite{chenruanorbcohom}).  Throughout, we
use real coefficients.

\end{definition}

Assume $Q$ is oriented, inducing an orientation of $\tilde{Q}_\Gamma$.  Bestow $\tilde{Q}_\Gamma$ with a Riemannian metric and metric
connection $\tilde{\omega}$ with curvature $\tilde{\Omega}$.  Let $E(\tilde{\Omega})$ denote
the Euler curvature form (see \cite{satake2} or \cite{seaton1}).  Note that if $(1)$ denotes the
$\approx$-class of the trivial homomorphism into any isotropy group, then as all such homomorphisms are clearly elements of the same $\approx$-class, $\tilde{Q}_{(1)}$ is clearly diffeomorphic to $Q$.  Restricting to the connected component $\tilde{Q}_{(1)}$, we have a metric
connection $\omega$ on $Q$ with curvature $\Omega$ and Euler curvature form $E(\Omega)$.

\begin{definition}[$\Gamma$-Euler-Satake class]
\label{def-eulerclass}

Let $Q$ be an oriented orbifold of dimension $n$ and $\Gamma$ a
finitely generated group. Let $e_{ES}(Q) \in H^n(Q)$ denote the
cohomology class of the Euler curvature form on $Q$.  We refer to it
as the {\bf Euler-Satake class of $Q$}. Let $e_\Gamma^{ES} (Q) \in
H_\Gamma^\ast (Q)$ denote the cohomology class of the Euler
curvature form on $\tilde{Q}_\Gamma$, called the {\bf
$\Gamma$-Euler-Satake class of $Q$}.

\end{definition}

Note that the Euler-Satake class and $\Gamma$-Euler-Satake class can
be defined in the obvious analogous manner for any good, oriented
orbifold vector bundle over $Q$.  Moreover, for bad, oriented
orbifold vector bundles, we can use the techniques in
\cite{seaton3}.  We are using the usual convention that the
$\Gamma$-Euler-Satake class of $Q$ indicates the
$\Gamma$-Euler-Satake class of the tangent bundle of $Q$.

The class $e_{ES}(Q)$ is the cohomology class represented by the
Gauss-Bonnet integrand in \cite{satake2}, while $e_\Z^{ES}(Q)$ is
the cohomology class represented by the Gauss-Bonnet integrand in
\cite[Theorem 3.2]{seaton1}. We have
\begin{equation}
\label{eq-sumofeulersatakeclasses}
    e_\Gamma^{ES}(Q) = \sum\limits_{(\phi) \in T_Q^\Gamma}
        e_{ES} \left(\tilde{Q}_{(\phi)} \right)
\end{equation}
so that $e_\Gamma^{ES}(Q)$ is generally not a homogeneous cohomology
class.  Satake's Gauss-Bonnet theorem for orbifolds implies that if
$Q$ is compact, then
\[
    \left\langle e_{ES}\left(\tilde{Q}_{(\phi)}\right);
    \left[\tilde{Q}_{(\phi)}\right]\right\rangle
    =
    \chi_{ES}\left(\tilde{Q}_{(\phi)}\right).
\]
By $\left\langle e_{ES}\left(\tilde{Q}_{(\phi)}\right);
\left[\tilde{Q}_{(\phi)}\right]\right\rangle$, we mean the integral
of any differential form representing
$e_{ES}\left(\tilde{Q}_{(\phi)}\right)$ on $\tilde{Q}_{(\phi)}$.
In particular, $e_\Gamma^{ES}(Q) = 0$ implies that $\chi_{ES}\left(\tilde{Q}_{(\phi)}\right) = 0$ for each $(\phi) \in T_Q^\Gamma$.

Conversely, since the top cohomology group of each $\Gamma$-sector
is isomorphic to $\R$ in the oriented case (see \cite[page 34 and
Theorem 2.13]{ademleidaruan}),
$\chi_{ES}\left(\tilde{Q}_{(\phi)}\right) = 0$ implies that
$e_{ES}\left(\tilde{Q}_{(\phi)}\right) = 0$.  If this is true for
each $(\phi)\in T_Q^\Gamma$, then by Equation
\ref{eq-sumofeulersatakeclasses}, $e_\Gamma^{ES}(Q) = 0$.  We
summarize this observation as follows.

\begin{lemma}
\label{lem-charsvanishiffclassvanishes}

Let $Q$ be a closed, oriented orbifold and $\Gamma$ a finitely generated group.  Then $e_\Gamma^{ES}(Q) = 0$ if and only if $\chi_{ES}\left(\tilde{Q}_{(\phi)}\right) = 0$ for each $(\phi) \in T_Q^\Gamma$.

\end{lemma}


\section{Topological Properties of the $\Gamma$-Sectors}
\label{sec-topologicalproperties}

As above, $Q$ is an $n$-dimensional orbifold whose orbifold
structure is given by the groupoid $\mathcal{G}$.  Recall from
Subsection \ref{subsec-connectedcomponents} that $\pi :
\tilde{Q}_\Gamma \rightarrow Q$ is the smooth map
$\pi(\mathcal{G}\phi_x) = \sigma(x)$. Note that for each $p \in Q$,
$\pi^{-1}(p)$ is finite, as $\mbox{HOM}(\Gamma, G_x)$ is finite for
each $x \in \sigma^{-1}(p)$ and the action of an $h \in G_1$ with
$t(h) = x$ identifies each element of $\mbox{HOM}(\Gamma, G_{t(x)})$
with an element of $\mbox{HOM}(\Gamma, G_x)$.

\begin{lemma}
\label{lem-gammasectorsclosed}

Suppose $Q$ is closed and $\Gamma$ is a finitely generated group.
Then each $\Gamma$-sector is a closed orbifold, and the image of each $\Gamma$-sector under $\pi$ is
a compact subset of $Q$.

\end{lemma}

\begin{proof}

Pick $(\phi) \in T_\Gamma^Q$.
We first claim that $\pi\left(\tilde{Q}_{(\phi)}\right)$ is a
compact subset of $Q$.  For each $p \in Q$, pick a linear chart at
$x \in G_0$, $\sigma(x) = p$,  with image $U_p \subseteq Q$.  For
each chart with domain $V_x$, let $V_x^\prime$ be a $G_x$-invariant
ball about $0 \in V_x$ such that $\overline{V_x^\prime} \subset
V_x$.  Let $U_p^\prime = \pi_p(V_x^\prime)$.  Then the $U_p^\prime$
form an open cover of $Q$.
As $Q$ is compact, there is a finite subcover $\{ U_{p_i}^\prime \}$
for $i = 1, 2, \ldots , k$ covered by linear charts at points $x_i$
with $\sigma(x_i) = p_i$. Then
\[
    \pi\left(\tilde{Q}_{(\phi)}\right) =
    \bigsqcup\limits_{\{i: \exists \phi_{x_i} \in (\phi)\}} \pi_{x_i}
    \left(\overline{V_{x_i}^\prime} \cap V_{x_i}^{\langle \phi_{x_i} \rangle} \right),
\]
which is a finite union of closed sets and hence closed and compact.

Now, suppose $(\phi_x)_i$ is a sequence in $\tilde{Q}_{(\phi)}$.
Then $\pi[(\phi_x)_i] = p_i$ is a sequence in the compact space
$\pi\left(\tilde{Q}_{(\phi)}\right)$, implying that it has a
subsequence $p_{i_j}$ with limit $p \in
\pi\left(\tilde{Q}_{(\phi)}\right)$.  As $\pi^{-1}(p) \cap
\tilde{Q}_{(\phi)}$ is a finite set, there is at least one $\phi_y
\in \pi^{-1}(p) \cap \tilde{Q}_{(\phi)}$ such that every
neighborhood of $\phi_y$ contains an infinite number of the
$(\phi_x)_{i_j}$.  It follows that there is a subsequence of
$(\phi_x)_{i_j}$ that converges to $\phi_y$, and that
$\tilde{Q}_{(\phi)}$ is compact.

\end{proof}

We define a relation on $T_Q^\Gamma$ as follows.  We say that
$(\psi) \leq (\phi)$ or equivalently $(\phi) \geq (\psi)$ if $\pi
\left(\tilde{Q}_{(\psi)} \right) \subseteq \pi \left(
\tilde{Q}_{(\phi)} \right)$, and $(\phi) \equiv (\psi)$ if $\pi
\left( \tilde{Q}_{(\phi)} \right) = \pi \left( \tilde{Q}_{(\psi)}
\right)$. It can happen that $(\phi) \equiv (\psi)$ with $(\phi)
\neq (\psi)$. However, we can consider $\leq$ a partial order of the
equivalence classes of $\Gamma$-sectors under the obvious
equivalence relation $\equiv$. By $(\psi) < (\phi)$, then, we will
mean that $(\psi) \leq (\phi)$ and $(\psi) \not\equiv (\phi)$.

The following technical lemma demonstrates that the partial order
$\leq$ can be understood completely locally.

\begin{lemma}
\label{lem-orderpreserved}

Let $Q$ be an orbifold and $\Gamma$ a finitely generated group.
Suppose there is an $x \in G_0$ and homomorphisms $\phi_x, \psi_x :
\Gamma \rightarrow G_x$ in $\approx$-classes $(\phi)$ and $(\psi)$,
respectively.  Let $\{ V_x, G_x, \pi_x \}$ be any linear chart for
$Q$ at $x$.
\newcounter{Ccount}
\begin{list}{\roman{Ccount}.}{\usecounter{Ccount}}
\item   If $V_x^{\langle \psi_x \rangle} \subseteq V_x^{\langle \phi_x
        \rangle}$, then $(\psi) \leq (\phi)$.
\item   If $V_x^{\langle \psi_x \rangle} \subset V_x^{\langle \phi_x
        \rangle}$, then $(\psi) < (\phi)$.
\item   If $V_x^{\langle \psi_x \rangle} = V_x^{\langle \phi_x \rangle}$,
        then $(\psi) \equiv (\phi)$.
\end{list}
\end{lemma}

\begin{proof}

First, we note that if any of the containment hypotheses involving
$V^{\langle \psi_x \rangle}$ and $V^{\langle \phi_x \rangle}$ are
true for any linear chart at $x$, then they are true for every
linear chart at $x$.  This follows from the fact that $G_x$ acts
linearly in every such chart so that these spaces are subspaces of
each $V_x \ni 0$; of course, subspaces are determined by their
intersection with any neighborhood of the origin.

If $h \in G_1$ with $s(h) = x$, then $h$ defines an equivalent
linear orbifold chart for $Q$ at $t(h)$ of the form $\left\{ t\circ
s_h^{-1}(V_x), hG_xh^{-1}, \pi_{t(h)} \} = \{ V_{t(h)}, G_{t(h)},
\pi_{t(h)} \right\}$. It is clear that
\begin{equation}
\label{eq-orderpres1}
    V_x^{\langle \psi_x \rangle} \subseteq V_x^{\langle \phi_x \rangle}
    \Longleftrightarrow
    V_{t(h)}^{\langle h\psi_x \rangle} \subseteq V_{t(h)}^{\langle h\phi_x \rangle}
\end{equation}
and
\begin{equation}
\label{eq-orderpres2}
    V_x^{\langle \psi_x \rangle} \subset V_x^{\langle \phi_x \rangle}
    \Longleftrightarrow
    V_{t(h)}^{\langle h\psi_x \rangle} \subset V_{t(h)}^{\langle h\phi_x
    \rangle}.
\end{equation}
We begin by showing that the containment hypotheses are preserved by a local covering in either direction.

Suppose there is a $\psi_y \in \mathcal{S}_\mathcal{G}^\Gamma$ with
$\psi_x \stackrel{loc}{\curvearrowright} \psi_y$, and then by Lemma
\ref{lem-localconglocalequal}, there is a linear chart $\{ V_x, G_x,
\pi_x \}$ at $x$ and an element of the orbit $\mathcal{G}\psi_y$ (which we assume, without loss of generality by Equations \ref{eq-orderpres1} and \ref{eq-orderpres2}, is equal to $\psi_y$) such that $y \in V_x$, and $\xi_x^y \circ \psi_y
= \psi_x$.  It follows, in particular, that $y \in V_x^{\langle
\psi_x \rangle}$.  Pick $\phi_x : \Gamma
\rightarrow G_x$ and assume that $V_x^{\langle \psi_x \rangle}
\subseteq V_x^{\langle \phi_x \rangle}$, and then as $y \in V_x^{\langle \phi_x \rangle}$,
$\xi_x^y(G_y)$ contains $\Imt \phi_x$ as a subgroup. Recalling that
$\xi_x^y$ is injective, define
\[
\begin{array}{rccl}
    \phi_y  :&     \Gamma  &\longrightarrow&       G_y         \\\\
            :&     \gamma  &\longmapsto&
                    (\xi_x^y)^{-1}[\phi_x(\gamma)].
\end{array}
\]
Then $\phi_y \in \mathcal{S}_\mathcal{G}^\Gamma$.  As $\xi_x^y \circ
\phi_y = \phi_x$, it is clear that $\phi_x
\stackrel{loc}{\curvearrowright} \phi_y$.  Pick a linear chart $\{
V_y, G_y, \pi_y \}$ at $y$, and assume by shrinking $V_y$ if
necessary that $V_y \subseteq V_x \subseteq G_0$.  Then it follows
from the construction that $V_y^{\langle \phi_y \rangle} = V_y \cap
V_x^{\langle \phi_x \rangle}$ and $V_y^{\langle \psi_y \rangle} =
V_y \cap V_x^{\langle \psi_x \rangle}$.  Therefore, $V_y^{\langle
\psi_y \rangle} \subseteq V_y^{\langle \phi_y \rangle}$.  Note
further that, if we assume a strict inclusion $V_x^{\langle \psi_x
\rangle} \subset V_x^{\langle \phi_x \rangle}$, then the as both
sets are subspaces of $V_x$, which is diffeomorphic to $\R^n$, and
as $V_y$ is an open subset of $V_x$, it follows that $V_y^{\langle
\psi_y \rangle} \subset V_y^{\langle \phi_y \rangle}$ is strict.

On the other hand, if there is a $\psi_y \in
\mathcal{S}_\mathcal{G}^\Gamma$ with $\psi_y
\stackrel{loc}{\curvearrowright} \psi_x$, then by Lemma
\ref{lem-localconglocalequal}, there is a linear chart $\{ V_y, G_y,
\pi_y \}$ at $y$ and a representative of $\mathcal{G}\psi_x$ (which
we assume, again without loss of generality, is equal to $\psi_x$)
such that $x \in V_y$, and $\xi_y^x \circ \psi_x = \psi_y$.  Let $\{
V_x, G_x, \pi_x \}$ be a linear chart at $x$; by shrinking $V_x$ if
necessary, we assume that $V_x \subseteq V_y \subseteq G_0$.  Assume
that $V_x^{\langle \psi_x \rangle} \subseteq V_x^{\langle \phi_x
\rangle}$ for some $\phi_x : \Gamma \rightarrow G_x$, and note that
as $\xi_y^x : G_x \rightarrow G_y$ extends the $G_x$-action on $V_x$
to all of $V_y$, $V_x^{\langle \psi_x \rangle} = V_x \cap
V_y^{\langle \psi_y \rangle}$.  Defining $\phi_x = \xi_y^x \circ
\phi_y$, we have that $\phi_x \in \mathcal{S}_\mathcal{G}^\Gamma$
and $\phi_y \stackrel{loc}{\curvearrowright} \phi_x$.  As above, we
have $V_x^{\langle \phi_x \rangle} = V_x \cap V_y^{\langle \phi_y
\rangle}$, so that $V_y^{\langle \psi_y \rangle} \subseteq
V_y^{\langle \phi_y\rangle}$.  Again, as the sets in question are
subspaces, if $V_x^{\langle \psi_x \rangle} \subset V_x^{\langle
\phi_x \rangle}$, then $V_y^{\langle \psi_y \rangle} \subset
V_y^{\langle \phi_y\rangle}$.

Now, suppose $V^{\langle \psi_x \rangle} \subseteq V^{\langle \phi_x
\rangle}$.  Let $q \in \pi\left(\tilde{Q}_{(\psi_x)}\right)$, and
then there is a $y \in G_0$ with $\sigma(y) = q$ and a $\psi_y \in
\mathcal{S}_\mathcal{G}^\Gamma$ with $\psi_y \approx \psi_x$.  By
the definition of $\approx$,  there is a finite sequence
$\psi_{x_0}, \psi_{x_1}, \ldots , \psi_{x_l}$ such that $\psi_{x_0}
\in {\mathcal G}\psi_x$, $\psi_{x_l} \in {\mathcal G}\psi_y$, and
for each $i$, $\psi_{x_i} \stackrel{loc}{\curvearrowright}
\psi_{x_{i+1}}$ or $\psi_{x_{i+1}} \stackrel{loc}{\curvearrowright}
\psi_{x_i}$.  By Equation \ref{eq-orderpres1}, we can assume that
$\psi_{x_0} = \psi_x$ and $\psi_{x_l} = \psi_y$.  Applying the above
arguments for each $i$, we have that there is a sequence
$\phi_{x_0}, \phi_{x_1}, \ldots, \phi_{x_l}$ such that $\phi_{x_0} =
\phi_x$, $\phi_{x_l} = \phi_y$, and for each $i$, $\phi_{x_i}
\stackrel{loc}{\curvearrowright} \phi_{x_{i+1}}$ or $\phi_{x_{i+1}}
\stackrel{loc}{\curvearrowright} \phi_{x_i}$.  At each step,
$V^{\langle \psi_{x_i} \rangle} \subseteq V^{\langle \phi_{x_i}
\rangle}$ implies that $V^{\langle \psi_{x_{i+1}} \rangle} \subseteq
V^{\langle \phi_{x_{i+1}} \rangle}$.  It follows that $q \in
\pi\left(\tilde{Q}_{(\phi)}\right)$, proving {\it (i)}.

To prove {\it (ii)}, we apply Equation \ref{eq-orderpres2} and note
that it was shown above that $V^{\langle \psi_{x_i} \rangle} \subset
V^{\langle \phi_{x_i} \rangle}$ implies that $V^{\langle
\psi_{x_{i+1}} \rangle} \subset V^{\langle \phi_{x_{i+1}} \rangle}$
for each $i$.
To prove {\it (iii)}, we simply apply {\it (i)} to $\psi_x$ and
$\phi_x$ and then reverse their roles.

\end{proof}

\begin{lemma}
\label{lem-opendensesubsetconnected}

Let $Q$ be an orbifold and $\Gamma$ a finitely generated group.
For each $(\phi) \in T_Q^\Gamma$, the set
\[
    \pi\left(\tilde{Q}_{(\phi)} \right) \backslash
    \bigcup\limits_{(\psi) < (\phi)}
    \pi\left(\tilde{Q}_{(\psi)}\right)
\]
is connected.

\end{lemma}

\begin{proof}

Pick an $x \in G_0$ and $\phi_x \in
\mathcal{S}_\mathcal{G}^\Gamma$ such that $\phi_x \in (\phi)$.  Then a
chart for $\tilde{Q}_{(\phi)}$ at $\phi_x$ is of the form $\left\{
V_x^{\langle \phi_x \rangle}, C_{G_x}(\phi_x), \pi_x^{\phi_x} \right\}$.
By \cite[Theorem 4.3.2, page 158]{pflaum}, there is an open and
dense subset $\mathcal{O}$ of $V_x^{\langle \phi_x \rangle}$ of principal orbit
type with respect to the $C_{G_x}(\phi_x)$-action such that
$\mathcal{O}/C_{G_x}(\phi_x)$ is connected.

If some point $y \in V_x^{\langle \phi_x \rangle}$ satisfies $\pi
\circ \pi_x^{\phi_x} (y) \in\pi\left(\tilde{Q}_{(\psi)}\right)$ for
some $(\psi) < (\phi)$, then there is a $\psi_y \in (\psi)$.  This
implies that the group $\langle \Imt \psi_y, \Imt \phi_y \rangle
\leq G_y$ where $\phi_y = (\xi_x^y)^{-1}\circ \phi_x$.  It follows
that $y$ is not an element of the principal orbit type; otherwise,
$\Imt(\xi_x^y \circ \psi_y)$ acts trivially on $V_x^{\langle \phi_x
\rangle}$, implying by Lemma \ref{lem-orderpreserved} that $(\psi)
\equiv (\phi)$ (contradicting the fact that $(\psi) < (\phi)$).
Conversely, if $y$ is not of principal orbit type, then picking
a surjective $\psi_y$ onto $G_y$ clearly defines a class $(\psi)$
with $(\psi) < (\phi)$. Hence, we have that in the image of each
local chart, $\tilde{Q}_{(\phi)}\backslash \bigcup\limits_{(\psi) <
(\phi)} \pi^{-1} \circ \pi\left(\tilde{Q}_{(\psi)}\right)$
corresponds to the connected set of points with principal
$C_{G_x}(\phi_x)$-orbit type. As $\tilde{Q}_{(\phi)}$ is connected,
so that any two points can be connected by a path covered by such
charts, this implies that $\tilde{Q}_{(\phi)}\backslash
\bigcup\limits_{(\psi) < (\phi)} \pi^{-1} \circ
\pi\left(\tilde{Q}_{(\psi)}\right)$ is connected. Note that if $q
\in \pi\left(\tilde{Q}_{(\psi)}\right)$ for some $(\psi) < (\phi)$,
then the isotropy groups of the points in $\pi^{-1}(q) \cap
\tilde{Q}_{(\phi)}$ are isomorphic.  Therefore, each such point is
contained in $\tilde{Q}_{(\psi^\prime)}$ for some $(\psi^\prime) <
(\phi)$, and
\[
    \pi\left(\tilde{Q}_{(\phi)}\backslash \bigcup\limits_{(\psi) < (\phi)}
    \pi^{-1} \circ \pi\left(\tilde{Q}_{(\psi)}\right)\right) =
    \pi\left(\tilde{Q}_{(\phi)} \right) \backslash
    \bigcup\limits_{(\psi) < (\phi)}
    \pi\left(\tilde{Q}_{(\psi)}\right)
\]
is the continuous image of a connected set, hence connected.

\end{proof}

We note the following, which is a trivial consequence of Lemma
\ref{lem-finitegammasectors}.

\begin{lemma}
\label{lem-existminimalsectors}

Let $Q$ be a closed orbifold and $\Gamma$ a finitely generated
group. For each $\Gamma$-sector $\tilde{Q}_{(\phi)}$ of $Q$, there
is a $\Gamma$-sector $\tilde{Q}_{(\psi)}$ of $Q$ such that $(\psi)
\leq (\phi)$ that represents a minimal $\equiv$-class with respect
to $\leq$. In other words, if $(\phi^\prime) \leq (\psi)$, then
$(\phi^\prime) \equiv (\psi)$.

\end{lemma}

We will abuse language slightly and say that $(\psi)$ is minimal
with respect to $\leq$.  By this, we mean that the $\equiv$-class of
$(\psi)$ is minimal.

\begin{definition}[Covering the local groups]

We say that the group $\Gamma$ {\bf covers the local groups of $Q$}
if, for each subgroup $H$ of each isotropy group $G_x$ of $Q$,
there is a homomorphism $\phi_x : \Gamma \rightarrow G_x$ with $\Imt
\phi_x = H$.

\end{definition}

We note that for every compact orbifold $Q$, there is a finite group
that covers the local groups of $Q$.  See the the proof of Lemma
\ref{lem-finitegammasectors}; for each $i = 1, 2, \ldots, k$, let
$\{ H_{i,j} : i = 1, 2, \ldots, l_i \}$ be a collection of all of
the nontrivial subgroups of $G_{x_i}$.  Then
\[
    \bigoplus\limits_{i=1}^k \bigoplus\limits_{j=1}^{l_i} H_{i,j}
\]
covers the local groups of $Q$.  Similarly, if $Q$ is any orbifold
such that the number of generators in a presentation of an isotropy
group of $Q$ is bounded by $d$, then the free group with $d$ generators $\mathbb{F}_d$ covers the local
groups of $Q$.

\begin{lemma}
\label{lem-intersectionssectors}

Suppose $\Gamma$ covers the local groups of $Q$ and $\pi \left(
\tilde{Q}_{(\phi)} \right) \cap \pi \left( \tilde{Q}_{(\psi)}
\right) \neq \emptyset$. Then there is a $(\psi^\prime) \in
T_Q^\Gamma$ with $(\psi^\prime) \leq (\phi)$, $(\psi^\prime) \leq
(\psi)$, and $\pi \left( \tilde{Q}_{(\psi^\prime)} \right) \subseteq
\pi \left( \tilde{Q}_{(\phi)} \right) \cap \pi \left(
\tilde{Q}_{(\psi)} \right)$.  Moreover, $\pi \left(
\tilde{Q}_{(\phi)} \right) \cap \pi \left( \tilde{Q}_{(\psi)}
\right)$ is a union of the image of such sectors; i.e.
\[
    \pi \left( \tilde{Q}_{(\phi)} \right) \cap \pi \left(
    \tilde{Q}_{(\psi)} \right)
    = \bigcup\limits_{(\psi^\prime) \leq (\phi), (\psi^\prime) \leq
    (\psi)} \pi\left(\tilde{Q}_{(\psi^\prime)}\right).
\]

\end{lemma}

Note that it is possible that $(\psi^\prime) \equiv (\phi)$,
$(\psi^\prime) \equiv (\psi)$, or both.

\begin{proof}

Pick $p \in \pi \left( \tilde{Q}_{(\phi)} \right) \cap \pi \left(
\tilde{Q}_{(\psi)} \right)$ and $x \in G_0$ with $\sigma(x) = p$.
Then there are $\phi_x, \psi_x \in \mathcal{S}_\mathcal{G}^\Gamma$
with $\phi_x \in (\phi)$, $\psi_x \in (\psi)$.  As $\Gamma$ covers
the local groups of $Q$, let $\psi_x^\prime : \Gamma \rightarrow
\langle \Imt \phi_x, \Imt \psi_x \rangle$ be surjective.  Then
clearly $V_x^{\langle \psi_x^\prime \rangle} \subseteq V_x^{\langle
\phi_x \rangle}$ and $V_x^{\langle \psi_x^\prime \rangle} \subseteq
V_x^{\langle \psi_x \rangle}$.  By Lemma \ref{lem-orderpreserved},
letting $(\psi^\prime)$ denote the $\approx$-class of
$\psi_x^\prime$ as usual, $(\psi^\prime) \leq (\phi)$ and
$(\psi^\prime) \leq (\psi)$.

This construction can be performed for each $p \in \pi \left(
\tilde{Q}_{(\phi)} \right) \cap \pi \left( \tilde{Q}_{(\psi)}
\right)$, so that any such $p$ is clearly contained in some
$\pi\left(\tilde{Q}_{(\psi^\prime)}\right)$ with $(\psi^\prime) \leq
(\phi)$, $(\psi^\prime) \leq (\psi)$.

\end{proof}

We note that in the case that $Q$ is an {\bf abelian orbifold}, i.e.
if each of the $G_x$ are abelian, then the restriction
$\pi|_{\tilde{Q}_{(\phi)}}$ of $\pi$ to any $\Gamma$-sector is
injective, and hence an embedding of $\tilde{Q}_{(\phi)}$ into $Q$
as a suborbifold. This follows from the fact that $C_{G(\phi_x)} =
G_x$ for each $\phi_x \in \mathcal{S}_\mathcal{G}^\Gamma$.  In
general, however, $\pi|_{\tilde{Q}_{(\phi)}}$ will fail to be
injective.  In the Lemmas \ref{lem-branchedcover} and
\ref{lem-minsectorsmanifolds}, we show that
$\pi|_{\tilde{Q}_{(\phi)}}$ is a sort of singular finite covering
space, and its singularities occur precisely on the images of
$\Gamma$-sectors $\tilde{Q}_{(\psi)}$ with $(\psi) < (\phi)$. When
there are no such sectors, $\pi|_{\tilde{Q}_{(\phi)}}$ is a covering
space of smooth manifolds.

\begin{lemma}
\label{lem-branchedcover}

Suppose $\Gamma$ covers the local groups of the orbifold $Q$.  Let
$\tilde{Q}_{(\phi)}$ have dimension $k$ and let $p \in \pi
\left(\tilde{Q}_{(\phi)} \right)$. One of the following is true.
\newcounter{Bcount}
\begin{list}{\roman{Bcount}.}{\usecounter{Bcount}}
\item   The point $p$ is contained in $\pi \left(\tilde{Q}_{(\psi)} \right)$ for some
        $(\psi) < (\phi)$.
\item   There is a neighborhood $W$ of $p$ in
        $\pi\left(\tilde{Q}_{(\phi)}\right)$ diffeomorphic to $\R^k$
        such that $\pi^{-1}(W) \cap \tilde{Q}_{(\phi)}$ is a finite number of disjoint sets
        diffeomorphic to $W$.
\end{list}
In particular, the set $\pi\left(\tilde{Q}_{(\phi)} \right)
\backslash \bigcup\limits_{(\psi) < (\phi)}
\pi\left(\tilde{Q}_{(\psi)}\right)$ is a smooth manifold equipped
with the trivial action of a finite group.

\end{lemma}

\begin{proof}

Pick $p \in \pi \left(\tilde{Q}_{(\phi)}\right)$ and $x \in G_0$
with $\sigma(x) = p$.  Let $\{ V_x, G_x, \pi_x \}$ be a linear chart
for $Q$ at $x$.  Then as $p = \sigma(x) \in
\left(\tilde{Q}_{(\phi)}\right)$, there is a $\phi_x \in
\mathcal{S}_\mathcal{G}^\Gamma$ that is a representative of
$(\phi)$.  By Lemma \ref{lem-orbifoldchartsforgammasectors}, a
linear chart at $\phi_x$ for the connected component
$\tilde{Q}_{(\phi)}$ of $\tilde{Q}_\Gamma$ is $\left\{ V_x^{\langle
\phi_x \rangle}, C_{G_x}(\phi_x), \pi_x^{\phi_x} \right\}$.  Note
that $V_x^{\langle \phi_x \rangle}$ is a subspace of $V_x$, and as
it forms an orbifold chart for $\tilde{Q}_{(\phi)}$, it has
dimension $k$.

Suppose $G_x$ does not act trivially on $V_x^{\langle \phi_x
\rangle}$ as a subset of $V_x$.  This means that there is a $g \in
G_x$ and a $y \in V_x^{\langle\phi_x\rangle}$ such that
$t[(\xi_x^y)^{-1}(g)] \neq y$ (of course, $t[(\xi_x^y)^{-1}(g)]$
need not be an element of $V_x^{\langle \phi_x \rangle}$).  As
$\Gamma$ covers the local groups of $Q$, let $\psi_x : \Gamma
\rightarrow G_x$ have image $\langle \Imt \phi_x, g \rangle$.  Then
$V_x^{\langle\psi_x\rangle}$ is a proper subspace of
$V_x^{\langle\phi_x\rangle}$ as it does not contain $y$.  By Lemma
\ref{lem-orderpreserved}, $(\psi_x) < (\phi_x)$, and ({\it i}) is
true.

Now, suppose $G_x$ acts trivially on $V_x^{\langle \phi_x \rangle}$.
Then $\pi_x^{\phi_x} : V_x^{\langle \phi_x \rangle} \rightarrow
\tilde{Q}_{(\phi)}$ is a diffeomorphism onto its image. Note that
$\pi \circ \pi_x^{\phi_x} = \pi_x$ on $V_x^{\langle \phi_x
\rangle}$, and $\pi_x$ is the quotient map by the trivial
$G_x$-action. Therefore, $\pi$ maps a neighborhood of $x$
diffeomorphic to $V_x^{\langle\phi_x\rangle}$ diffeomorphically onto
a neighborhood of $p$ in $\pi\left(\tilde{Q}_{(\phi)}\right)$.  Let
$W = \pi_x \left(V_x^{\langle\phi_x\rangle}\right)$ be this neighborhood. Any
other element of $\tilde{Q}_{(\phi)}$ in $\pi^{-1}(p)$ is of the
form $h\phi_x$ for an $h \in G_1$ with $s(h) = x$.  Then $\left\{
V_{t(x)}^{\langle h \phi_x \rangle}, C_{G_{t(h)}}(h\phi_x),
\pi_{t(h)}^{h\phi_x} \right\}$ is an equivalent orbifold chart for
$\tilde{Q}_{(\phi)}$ at $h\phi_x$.  Suppose $t(h) = x$ and
$V_{t(x)}^{\langle h \phi_x \rangle} \neq V_x^{\langle \phi_x
\rangle}$.  As $\Gamma$ covers the local groups of $Q$, let $\psi_x
: \Gamma \rightarrow G_x$ with image $\langle \Imt \phi_x, \Imt
h\phi_x \rangle$.  Then $V_x^{\langle \psi_x \rangle} =
V_{t(x)}^{\langle h \phi_x \rangle} \cap V_x^{\langle \phi_x
\rangle} \subset V_x^{\langle \phi_x \rangle}$ so that, by Lemma
\ref{lem-orderpreserved}, $(\psi) < (\phi)$.  Hence, if $t(h) = x$,
then either ({\it i}) is true, or $V_{t(x)}^{\langle h \phi_x
\rangle} = V_x^{\langle \phi_x \rangle}$.

If $t(h) \neq x$, then as the restriction $\mathcal{G}|_{V_x}$ is
isomorphic to $G_x \ltimes V_x$, $t(h) \not\in V_x$.  By shrinking
$V_x$ (and hence $W$) if necessary, we may assume that as subsets of
$G_0$, $V_x^{\langle\phi_x\rangle} \cap V_{t(h)}^{\langle
h\phi_x\rangle} = \emptyset$. Therefore,
\[
    \pi^{-1}(W) \cap \tilde{Q}_{(\phi)}
    =
    \bigsqcup\limits_{h \in s^{-1}(x),t(h)\neq x} \pi_{t(h)}^{h\phi_x}\left(V_{t(h)}^{\langle
    h \phi_x \rangle}\right).
\]
Note that as $\pi^{-1}(p)$ and hence $\pi^{-1}(p) \cap
\tilde{Q}_{(\phi)}$ is finite, the $\left\{ V_{t(x)}^{\langle h
\phi_x \rangle}, C_{G_{t(h)}}(h\phi_x), \pi_{t(h)}^{h\phi_x}
\right\}$ form linear charts for a finite number of open subsets of
$\tilde{Q}_{(\phi)}$. By the argument above, two sets
$V_{t(h_1)}^{\langle h_1\phi_x \rangle}$ and $V_{t(h_2)}^{\langle
h_2\phi_x \rangle}$ either are disjoint or coincide and yield
equivalent linear charts.  As each $V_{t(h)}^{\langle h \phi_x
\rangle}$ is a $k$-dimensional subspace of $V_{t(h)}$, and as
$\pi_{t(h)}^{h\phi_x}$ is simply the quotient map of the trivial
$C_{G_{t(h)}}(h\phi_x)$-action, ({\it ii}) is true.

Now, for every point $p \in \pi\left(\tilde{Q}_{(\phi)} \right)
\backslash \bigcup\limits_{(\psi) < (\phi)}
\pi\left(\tilde{Q}_{(\psi)}\right)$ there is a neighborhood $W_p$ of
$p$ in $\pi\left(\tilde{Q}_{(\phi)}\right)\backslash
\bigcup\limits_{(\psi) < (\phi)}
\pi\left(\tilde{Q}_{(\psi)}\right)$,  diffeomorphic to $\R^k$ such
that $\pi^{-1}(W_p)\cap \tilde{Q}_{(\phi)}$ is a finite number of
sets diffeomorphic to $W_p$. For each $\mathcal{G}\phi_x \in
\pi^{-1}(p) \cap \tilde{Q}_{(\phi)}$, there is a chart $\left\{
V_x^{\langle \phi_x \rangle}, C_{G_x}(\phi_x), \pi_x^{\phi_x}
\right\}$ for $\tilde{Q}_{(\phi)}$ at $\mathcal{G}\phi_x$ in which
the preimage of $W_p$ corresponds to $V_x^{\langle\phi_x\rangle}$
with trivial $C_{G_x}(\phi_x)$-action. As $\tilde{Q}_{(\phi)}$ is
connected, and as isotropy groups of objects in the same
$\mathcal{G}$-orbit are isomorphic, it is easy to see that the
isotropy group of each point in $\pi\left(\tilde{Q}_{(\phi)} \right)
\backslash \bigcup\limits_{(\psi) < (\phi)}
\pi\left(\tilde{Q}_{(\psi)}\right)$ is isomorphic. Moreover, each of
these manifold charts is the restriction of a linear orbifold chart
to an invariant subspace.  Hence, they patch together to give
$\pi\left(\tilde{Q}_{(\phi)} \right) \backslash
\bigcup\limits_{(\psi) < (\phi)} \pi\left(\tilde{Q}_{(\psi)}\right)$
the structure of a smooth manifold.

It follows that $\pi\left(\tilde{Q}_{(\phi)} \right) \backslash
\bigcup\limits_{(\psi) < (\phi)} \pi\left(\tilde{Q}_{(\psi)}\right)$
is a smooth manifold equipped with the trivial action of a finite
group; moreover, selecting $p \in \pi\left(\tilde{Q}_{(\phi)}
\right) \backslash \bigcup\limits_{(\psi) < (\phi)}
\pi\left(\tilde{Q}_{(\psi)}\right)$ and $x \in G_0$ such that $\sigma(x) = p$, that finite group is given by
$G_x$.

\end{proof}

\begin{lemma}
\label{lem-minsectorsmanifolds}

Suppose $\Gamma$ covers the local groups of $Q$.  If $(\phi)$ is a
minimal element of $T_Q^\Gamma$, then both $\tilde{Q}_{(\phi)}$ and
$\pi\left(\tilde{Q}_{(\phi)}\right)$ are smooth manifolds equipped
with the trivial action of a finite group.

\end{lemma}

\begin{proof}

That $\pi\left(\tilde{Q}_{(\phi)}\right)$ is a manifold follows
directly from Lemma \ref{lem-branchedcover}. Clearly in this case,
\[
    \pi\left(\tilde{Q}_{(\phi)} \right) \backslash
    \bigcup\limits_{(\psi) < (\phi)} \pi\left(\tilde{Q}_{(\psi)}\right)
    = \pi\left(\tilde{Q}_{(\phi)} \right).
\]

That $\tilde{Q}_{(\phi)}$ is a manifold follows from the fact that
the groups in the orbifold charts $\left\{ V_x^{\langle \phi_x
\rangle}, C_{G(\phi_x)}, \pi_x^{\phi_x} \right\}$ for
$\tilde{Q}_{(\phi)}$ act trivially.  Hence, $\tilde{Q}_{(\phi)}$ is
an orbifold in which every element of the local group acts
trivially, and hence the associated reduced orbifold is a smooth
manifold.

\end{proof}

Note in particular that if $(\phi)$ is minimal, then by Lemma
\ref{lem-branchedcover}, $\pi|_{\tilde{Q}_{(\phi)}}$ is a covering
map for its image.


\section{The Euler-Satake Characteristics and Classes as Complete Obstructions}
\label{sec-obstruction}

In this section, we use the constructions developed above to give a
necessary and sufficient condition for an orbifold to admit a
nonvanishing, smooth vector field.  Our main result is Theorem
\ref{thrm-closedobstruction} which deals with the case of a closed
orbifold; this is proven in Subsection \ref{subsec-closedcase}.  In Subsection \ref{subsec-boundarycase}, we prove Theorem \ref{thrm-boundarycase}, dealing
with the case of an orbifold with boundary, and Theorem
\ref{thrm-opencase}, dealing with a certain class of open orbifolds.

We start with two lemmas dealing with continuously extending and
smoothly approximating vector fields on closed orbifolds.

\begin{lemma}
\label{lem-extendvfields}

Let $Q$ be a closed orbifold and let $S \subseteq Q$ be closed.  A
continuous vector field on $S$ can be extended to a continuous
vector field on $Q$.

\end{lemma}

\begin{proof}

Suppose $X_0$ is a vector field defined on $S$.  Let $\{ (U_i, f_i)
\mid i = 1, 2, \ldots , m \}$ be a finite partition of unity for $Q$
composed of uniformized sets, each uniformized by a linear chart $\{ V_i,
G_i, \pi_i \}$.  For each $i$ such that $S \cap U_i$ is not empty,
we have that $\pi_i^{-1}(S)$ is a closed subset of $V_i$ and
$\pi_i^\ast X_0$ is a $G_i$-invariant vector field on $\pi_i^{-1}(S)$.
As $V_i$ is an open subset of $\R^n$, we can treat a vector field on
$\pi_i^{-1}(S) \subseteq V_i$ as $n$ $\R^n$-valued functions
$(\pi^\ast X_0)^j$, $j = 1, 2, \ldots , n$.  Extending each of these
functions to all of $V_i$ by the Tietsze Extension Theorem, we form
a vector field $Y_i$ on $V_i$ that extends $\pi_i^\ast X_0$.  Let
$X_i$ be the average of $Y_i$ over the $G_i$-action, i.e.
\[
    X_i = \frac{1}{|G_i|} \sum\limits_{g \in G_i} gY_i,
\]
and then as
$\pi_i^\ast X_0$ is $G_i$-invariant, $X_i$ also extends $\pi_i^\ast
X_0$.

For each $i$ such that $S \cap U_i = \emptyset$, let $Y_i$ be an
arbitrary vector field on $V_i$ and let $X_i$ be its average over
the $G_i$-action.

Since each $X_i$ is a $G_i$-invariant vector field on $V_i$, it
defines a vector field on $U_i$ (also denoted $X_i$).  The vector
field $X(p) = \sum\limits_{i = 1}^m f_i(p) X_i(p)$ is a continuous vector
field on $Q$ that extends $X_0$.

\end{proof}

\begin{lemma}
\label{lem-smoothvfldapprox}

Let $Q$ be a closed orbifold that admits a continuous vector field
that is nonvanishing on the closed set $S$. Then $Q$ admits a
${\mathcal C^\infty}$ vector field that is nonvanishing on the
closed set $S$.

\end{lemma}

\begin{proof}

Let $Y$ be a continuous vector field on $Q$ that is nonvanishing on
$S$. Fix a metric on $Q$ and let $M$ be the minimum value of $\|
Y(p) \|$ on $S$. Then $M > 0$.

For each $p \in Q$, pick a linear chart $\{ V_x, G_x, \pi_x \}$ at
some $x$ with $\sigma(x) = p$.  Then $\pi_x^\ast Y$ is a continuous
vector field so that there is a $G_x$-invariant open ball $W_x$
about $0 \in V_x$ such that
\[
    \| \pi_x^\ast Y(y) - \pi_x^\ast Y(0) \| < \frac{M}{2} .
\]
Note that $\pi_x^\ast Y(0)$ is a $G_x$-invariant vector field.  The
collection of the $\pi_x(W_x)$ form an open cover of $Q$, so let $\{
\pi_{x_i}(W_{x_i}) : i = 1, 2, \ldots , k \}$ be a finite subcover
with $\sigma(x_i) = p_i$ for each $i$. Let $\{ \rho_i : i = 1, 2,
\ldots , k \}$ be a partition of unity subordinate to this subcover
and define
\[
    X(p)   =   \sum\limits_{i = 1}^k \rho_i(p)Y(p_i) .
\]
Note that $X(p)$ is a smooth vector field on $Q$.  Then we have for
each $p \in Q$ that
\[
\begin{array}{rcl}
    \| Y(p) - X(p) \|
    &=&
    \left\| \sum\limits_{i = 1}^k \rho_i(p)Y(p) - \sum\limits_{i = 1}^k \rho_i(p)Y(p_i) \right\|
                                                                \\\\
    &=&
    \left\| \sum\limits_{\{i : p \in \mbox{\smaller supp}\: \rho_i \}} \rho_i(p)\left(Y(p) - Y(p_i)\right) \right\|
                                                                \\\\
    &\leq&
    \sum\limits_{\{i : p \in \mbox{\smaller supp}\: \rho_i \}} \rho_i(p)\left\| Y(p) - Y(p_i) \right\|
                                                                \\\\
    &<&
    \sum\limits_{\{ i : p \in \mbox{\smaller supp}\: \rho_i \}} \rho_i(p)\frac{M}{2}
                                                                \\\\
    &=&
    \frac{M}{2}.
\end{array}
\]
Therefore,
\[
\begin{array}{rcl}
    \| X(p) \|
    &=&
    \| Y(p) - (Y(p) - X(p)) \|
    \\\\
    &\geq& \| Y(p) \| -
    \left\| Y(p) - X(p) \right\|
    \\\\
    &>&     \| Y(p) \| - \frac{M}{2}
    \\\\
    &>&  \frac{M}{2} \:  > \: 0.
\end{array}
\]
Hence, $X$ is nonvanishing on $S$.

\end{proof}


\subsection{Closed Orbifolds}
\label{subsec-closedcase}

We turn to the proof of Theorem \ref{thrm-closedobstruction}.  One
direction of the theorem is true for any finitely generated group
$\Gamma$, so we state and prove it as Lemma
\ref{lem-nonvanimplieseulerclassvanishes}.  To prove the other
direction, we need to construct a smooth, nonvanishing vector field
$X$ on a closed orbifold $Q$ assuming that $\Gamma$ covers the local groups of $Q$, and $\chi_{ES}\left(\tilde{Q}_{(\phi)}\right) = 0$ for
each $\Gamma$-sector.  We will
construct $X$ on the $\Gamma$-sectors of $Q$ inductively using the
partial order $\leq$. To simplify the exposition, we will organize
this construction into two claims; Claim \ref{clm-closedbasecase} is
a base case and Claim \ref{clm-closedinductionstep} is the inductive
step.  The actual induction will be explained in the proof of the
theorem.

\begin{lemma}
\label{lem-nonvanimplieseulerclassvanishes}

Let $Q$ be a closed orbifold and $\Gamma$ a finitely generated group.  If $Q$ admits a nonvanishing, smooth
vector field $X$, then $\chi_{ES}\left(\tilde{Q}_{(\phi)}\right) = 0$ for
each $(\phi) \in T_Q^\Gamma$.

\end{lemma}

\begin{proof}

Suppose $Q$ admits a nonvanishing, smooth vector field $X$ and let
$\tilde{X}_\Gamma$ denote the extension of $X$ to $\tilde{Q}_\Gamma$
defined in Lemma \ref{lem-bundlesovergammasectors}.  Then
$\tilde{X}_\Gamma$ is smooth and nonvanishing.  Pick $(\phi) \in
T_Q^\Gamma$.  Letting
$\mbox{ind}^{orb}\left(\tilde{X}_\Gamma|_{\tilde{Q}_{(\phi)}};
\tilde{Q}_{(\phi)}\right)$ denote the index of the vector field
$\tilde{X}_\Gamma|_{\tilde{Q}_{(\phi)}}$ on the set
$\tilde{Q}_{(\phi)}$ in the orbifold sense (see \cite{satake2}), we
have
\[
    \mbox{ind}^{orb}\left(\tilde{X}_\Gamma|_{\tilde{Q}_{(\phi)}};
    \tilde{Q}_{(\phi)}\right) = 0.
\]
By the Poincar\'{e}-Hopf Theorem for closed
orbifolds in \cite{satake2}, then,
\[
\begin{array}{rcl}
    0
    &=&
    \mbox{ind}^{orb}\left(\tilde{X}_\Gamma|_{\tilde{Q}_{(\phi)}}; \tilde{Q}_{(\phi)}\right)
                            \\\\
    &=&
    \chi_{ES}\left(\tilde{Q}_{(\phi)}\right).
\end{array}
\]

\end{proof}

Now, we assume that $\Gamma$ covers
the local groups of the closed orbifold $Q$.

\begin{claim}[Base Case]
\label{clm-closedbasecase}

Let $Q$ be a closed orbifold and $\Gamma$ a finitely generated group
that covers the local groups of $Q$. If $\chi_{ES}\left(\tilde{Q}_{(\phi)}\right) = 0$ for each minimal $(\phi) \in T_Q^\Gamma$, then
there is a smooth vector field $X$ on $Q$ whose restriction to
$\pi\left(\tilde{Q}_{(\phi)}\right)$ for each minimal
$\Gamma$-sector $\tilde{Q}_{(\phi)}$ is nonvanishing.

\end{claim}

\begin{proof}

Let $(\phi)$ be a minimal element of $T_Q^\Gamma$.  Then
$\tilde{Q}_{(\phi)}$ and $\pi\left(\tilde{Q}_{(\phi)}\right)$ are
smooth manifolds equipped with the trivial action of a finite group
by Lemma \ref{lem-minsectorsmanifolds}.  As all isotropy groups of $\tilde{Q}_{(\phi)}$ are isomorphic so that the Euler-Satake
characteristic of $\tilde{Q}_{(\phi)}$ is simply its Euler
characteristic divided by the order of any isotropy group,
we see that $\chi_{top}\left(\tilde{Q}_{(\phi)}\right) = 0$.
Moreover, as $\pi|_{\tilde{Q}_{(\phi)}}$ is a covering map onto
$\pi\left(\tilde{Q}_{(\phi)}\right)$, it follows that
$\chi_{top}\left(\pi\left(\tilde{Q}_{(\phi)}\right)\right) = 0$.
Hence, it admits a smooth, nonvanishing vector field.

Noting that the images of minimal $\Gamma$-sectors are either
disjoint or coincide by Lemma \ref{lem-intersectionssectors}, we can
use this technique to construct a nonvanishing vector field on the
image of each minimal $\Gamma$-sector in $Q$. By Lemma
\ref{lem-gammasectorsclosed}, the union of the images of the minimal
sectors in $Q$ is a finite union of closed sets and hence closed.  Therefore, by Lemma
\ref{lem-extendvfields} we can extend to a vector field on all of $Q$ which, by Lemma \ref{lem-smoothvfldapprox} we may assume is smooth.

\end{proof}

\begin{claim}[Induction Step]
\label{clm-closedinductionstep}

Let $Q$ be a closed orbifold, $\Gamma$ a finitely generated group
that covers the local groups of $Q$, and assume $\chi_{ES}\left(\tilde{Q}_{(\phi)}\right) = 0$ for each $(\phi) \in T_Q^\Gamma$.  Let $(\phi) \in T_Q^\Gamma$, and suppose there is a continuous
vector field $X$ on $Q$ that restricts to a nonvanishing vector
field on $\bigsqcup\limits_{(\psi) < (\phi)}
\pi\left(\tilde{Q}_{(\psi)}\right)$.  Then there is a continuous
vector field $Y$ on $Q$ that does not vanish on
$\pi\left(\tilde{Q}_{(\phi)}\right)$ and coincides with $X$ on each
$\pi\left(\tilde{Q}_{(\psi)}\right)$ such that $(\psi) < (\phi)$.

\end{claim}

\begin{proof}

The zeros of $X|_{\pi \left(\tilde{Q}_{(\phi)}\right)}$ are
contained in the open set $\pi \left(\tilde{Q}_{(\phi)}\right)
\backslash \bigcup\limits_{(\psi) < (\phi)}
\pi\left(\tilde{Q}_{(\psi)}\right)$. By Lemmas
\ref{lem-opendensesubsetconnected} and \ref{lem-branchedcover}, this
set is a connected manifold.  Fix a point
\linebreak
$p \in \pi \left(\tilde{Q}_{(\phi)}\right) \backslash \bigcup\limits_{(\psi) <
(\phi)} \pi \left(\tilde{Q}_{(\psi)}\right)$ and an open
neighborhood $W$ of $p$ as in Lemma \ref{lem-branchedcover}.  We may
assume, by shrinking $W$ if necessary, that $W$ is contained in the
open set $\pi \left(\tilde{Q}_{(\phi)}\right) \backslash
\bigcup\limits_{(\psi) < (\phi)} \pi
\left(\tilde{Q}_{(\psi)}\right)$.  As $X$ does not vanish on the
closed set $\bigsqcup\limits_{(\psi) < (\phi)}
\pi\left(\tilde{Q}_{(\psi)}\right)$, we may continuously perturb
$X|_{\pi \left(\tilde{Q}_{(\phi)}\right)}$ away from each of the
$\pi\left(\tilde{Q}_{(\psi)}\right)$, so that we can assume the
zeros of $X$ are isolated and contained in the interior of a compact
set $K \subset W$.

We have that, $\pi^{-1}(K) \cap \tilde{Q}_{(\phi)}$ is a finite disjoint
union of sets diffeomorphic to $K$; say $\pi^{-1}(K)\cap
\tilde{Q}_{(\phi)} = \bigsqcup\limits_{j=1}^l J_i$ where each $J_i$
is diffeomorphic to $K$. Hence, $\tilde{X}_\Gamma$ restricts to a
continuous vector field $\tilde{X}_\Gamma|_{\tilde{Q}_{(\phi)}}$ on
$\tilde{Q}_{(\phi)}$ with only isolated zeros contained in the $J_i$
such that each of the $\tilde{X}_\Gamma|_{J_i}$ coincide with $X|_{K}$ via the
diffeomorphism between each $J_i$ and $K$. We have by the Poincar\'{e}-Hopf Theorem for closed orbifolds in \cite{satake2} that
\[
\begin{array}{rcl}
        0
        &=&
        \chi_{ES}\left(\tilde{Q}_{(\phi)}\right)                    \\\\
        &=&     \mbox{ind}^{orb}(\tilde{X}_\Gamma; \tilde{Q}_{(\phi)}) \\\\
        &=&     \sum\limits_{i=1}^l \mbox{ind}^{orb} (\tilde{X}_\Gamma; J_i)
                                                        \\\\
        &=&     l \left(\mbox{ind}^{orb}(X|_{\pi \left(\tilde{Q}_{(\phi)}\right)}; K)\right),
\end{array}
\]
so that $\mbox{ind}^{orb}(X|_{\pi \left(\tilde{Q}_{(\phi)}\right)}; K) = 0$.  By techniques
in \cite{gp}, $X|_{\pi \left(\tilde{Q}_{(\phi)}\right)}$ can be
perturbed continuously on an open set whose closure is contained in
$K$ resulting in a continuous, nonvanishing vector field on
$\pi\left(\tilde{Q}_{(\phi)}\right)$.  Applying Lemma
\ref{lem-extendvfields} with $S =
\pi\left(\tilde{Q}_{(\phi)}\right)$, we can extend to a continuous vector field
$Y$ on $Q$ that does not vanish on
$\pi\left(\tilde{Q}_{(\phi)}\right)$ and coincides with $X$ on each
$\pi\left(\tilde{Q}_{(\psi)}\right)$ with $(\psi) < (\phi)$.

\end{proof}

\begin{proof}[Proof of Theorem \ref{thrm-closedobstruction}]

Let $Q$ be a closed orbifold and $\Gamma$ a finitely generated group
that covers the local groups of $Q$.  Note that if $Q$ is oriented, then
$\chi_{ES}\left(\tilde{Q}_{(\phi)}\right) = 0$ for each $(\phi) \in T_Q^\Gamma$ is equivalent to $e_\Gamma^{ES}(Q) = 0$ by Lemma \ref{lem-charsvanishiffclassvanishes}.
Let $T_{\Gamma, 0}^Q \subseteq
T_Q^\Gamma$ denote the minimal elements. Let $T_{\Gamma, 1}^Q$
denote the set of $(\phi) \in T_Q^\Gamma \backslash T_{\Gamma, 0}^Q$
such that whenever $(\psi) < (\phi)$, $(\psi) \in T_{\Gamma, 0}^Q$.
Similarly, for each natural $j$, let $T_{\Gamma, j}^Q$ denote the
set of $(\phi) \in T_Q^\Gamma \backslash \bigcup\limits_{i = 0}^{j -
1} T_{\Gamma, i}^Q$ such that whenever $(\psi) < (\phi)$, $(\psi)
\in T_{\Gamma, i}^Q$ for some $i < j$.  Then by Lemma
\ref{lem-finitegammasectors}, there is an $m$ such that
\[
    T_Q^\Gamma = \bigcup\limits_{i=1}^m T_{\Gamma, i}^Q.
\]
In particular, $(1) \in T_{\Gamma, m}^Q$ where $(1)$ denotes the
$\approx$-class of the trivial homomorphism into any isotropy group.

By Claim \ref{clm-closedbasecase}, there is a smooth vector field
$X_0$ on $Q$ that does not vanish on $\bigsqcup\limits_{(\phi) \in
T_{\Gamma, 0}^Q} \pi\left(\tilde{Q}_{(\phi)}\right)$.  Pick $j$ with
$1 \leq j < m$ and assume that there is a continuous vector field
$Y_j$ on $Q$ that is nonvanishing on
$\pi\left(\tilde{Q}_{(\psi)}\right)$ for each $(\psi) \in T_{\Gamma,
j}^Q$.  This implies that $Y_j$ is nonvanishing on
$\pi\left(\tilde{Q}_{(\psi)}\right)$ for each $(\psi) \in T_{\Gamma,
i}^Q$ with $i \leq j$.

For each $(\phi) \in T_{\Gamma, j+1}^Q$, by Claim
\ref{clm-closedinductionstep} we can construct a continuous vector
field $Y_{(\phi)}$ on $Q$ such that the restriction
$Y_{(\phi)}|_{\pi\left(\tilde{Q}_{(\phi)}\right)}$ to
$\pi\left(\tilde{Q}_{(\phi)}\right)$ is nonvanishing.  If $(\phi),
(\phi^\prime) \in T_{\Gamma, j+1}^Q$, then $Y_{(\phi)}$ and
$Y_{(\phi^\prime)}$ need not coincide.  However, since they both
extend $Y_j$, it is clear that they coincide on
$\pi\left(\tilde{Q}_{(\psi)}\right)$ for each $(\psi) \in T_{\Gamma,
i}^Q$ with $i \leq j$.  Moreover, by Lemma
\ref{lem-intersectionssectors}, $\pi\left(\tilde{Q}_{(\phi)}\right)
\cap \pi\left(\tilde{Q}_{(\phi^\prime)}\right)$ is a union of such
$\pi\left(\tilde{Q}_{(\psi)}\right)$, so that $Y_{(\phi)}$ and
$Y_{(\phi^\prime)}$ coincide on $\pi\left(\tilde{Q}_{(\phi)}\right)
\cap \pi\left(\tilde{Q}_{(\phi^\prime)}\right)$. Hence, if we set
\[
    Y_{j+1}(p)  =   Y_{(\phi)}(p) \;\;\; \forall p \in \pi\left(\tilde{Q}_{(\phi)}
                \right), (\phi) \in T_{\Gamma, j+1}^Q,
\]
then $Y_{j+1}$ is a well-defined, continuous, nonvanishing vector
field on $\bigcup\limits_{(\phi) \in T_{\Gamma, j+1}^Q}
\pi\left(\tilde{Q}_{(\phi)} \right)$.  As this set is closed, we
apply Lemma \ref{lem-extendvfields} to extend $Y_{j+1}$ to a
continuous vector field (also denoted $Y_{j+1}$) on $Q$ that is
nonvanishing on $\pi\left(\tilde{Q}_{(\phi)}\right)$ for each
$(\phi) \in T_{\Gamma, j+1}^Q$

By induction, then, there is a continuous, nonvanishing vector field
$Y_m$ on $\tilde{Q}_{(1)}$, which is diffeomorphic to $Q$.  By Lemma
\ref{lem-smoothvfldapprox}, we can approximate $Y_m$ with a smooth,
nonvanishing vector field, completing the proof of Theorem
\ref{thrm-closedobstruction}.

\end{proof}

We end this subsection with an example of a closed orbifold that does not admit a nonvanishing vector field.
In this case, the obstruction is not detected when $\Gamma = \Z$ yet is detected for other choices of $\Gamma$.

\begin{example}
\label{ex-dihedralonS5}

Let $\R^6$ have standard basis $\{ e_1, e_2, e_3, e_4, e_5, e_6 \}$, and let the dihedral group $D_6$ act on the sphere $S^5 \subset \R^6$ as follows.  We let $a$ denote the permutation $(123)$ acting on the basis elements and
\[
    b   =   \left[ \begin{array}{cccccc}
                0   &   1   &   0   &   0   &   0   &   0       \\
                1   &   0   &   0   &   0   &   0   &   0       \\
                0   &   0   &   1   &   0   &   0   &   0       \\
                0   &   0   &   0   &   -1  &   0   &   0       \\
                0   &   0   &   0   &   0   &   1   &   0       \\
                0   &   0   &   0   &   0   &   0   &   1
                \end{array} \right]
\]
act by a permutation along with multiplying the fourth basis element by $-1$.  One checks that $\langle a, b\rangle$ is isomorphic to the dihedral group $D_6$.  We let $Q = S^5/D_6$, and then the orbifold groupoid $\mathcal{G} = S^5 \ltimes D_6$ is a representative of the orbifold structure of $Q$.

Letting $\Gamma = \Z$ with generator $\gamma$, there are three $\Gamma$-sectors.  The first corresponds to the homomorphism $\gamma \mapsto 1$ at each point, and is clearly diffeomorphic to $Q$.  The second corresponds to $\gamma \mapsto a$ over points stabilized by $a$, and is given by $S^3 = \mbox{Span}\;\{ e_1 + e_2 + e_3, e_4, e_5, e_6\} \cap S^5$ with trivial $\Z^3$-action.  The third corresponds to $\gamma \mapsto b$ over points stabilized by $b$, and is given by $S^3 = \mbox{Span}\;\{ e_1 + e_2, e_3, e_5, e_6\} \cap S^5$ with trivial $\Z^2$-action.  We note that the $\Z$-sectors correspond to the inertia orbifold; the Euler-Satake characteristics of each of these sectors is equal to zero, as they are all odd-dimensional orbifolds (see \cite[Theorem 4]{satake2}).

Now, let $\Gamma = \F_2$, the free group with generators $\gamma_1$ and $\gamma_2$.  Designating a homomorphism $\F_2 \rightarrow D_6$ by $(g_1, g_2)$ where $\gamma_i \mapsto g_i$, we have the following conjugacy classes:
\[
\begin{array}{c}
    \{(1,1) \}
\end{array}
\]
mapping into the isotropy group over every point,
\[
\begin{array}{c}
    \{ (1, a), (1, a^2) \};         \\
    \{ (a, 1), (a^2, 1) \}; \\
    \{ (a, a), (a^2, a^2) \}; \\
    \{ (a, a^2), (a^2, a) \}
\end{array}
\]
mapping into the isotropy group over every point in $S^3 = \mbox{Span}\;\{ e_1 + e_2 + e_3, e_4, e_5, e_6\} \cap S^5$,
\[
\begin{array}{c}
    \{ (1, b), (1, ab), (1, a^2b) \}; \\
    \{ (b, 1), (ab, 1), (a^2b, 1) \}; \\
    \{ (b, b), (ab, ab), (a^2b, a^2b) \}
\end{array}
\]
mapping into the isotropy group over every point in $S^3 = \mbox{Span}\;\{ e_1 + e_2, e_3, e_5, e_6\} \cap S^5$ (or a representative of the orbit of this set), and
\[
\begin{array}{c}
    \{ (b, ab), (ab, a^2b), (a^2b, b), (ab, b), (a^2b, ab), (b, a^2b) \}; \\
    \{ (a, b), (a, ab), (a, a^2b), (a^2, b), (a^2, ab), (a^2, a^2b) \}; \\
    \{ (b, a), (ab, a), (a^2b, a), (b, a^2), (ab, a^2), (a^2b, a^2) \};
\end{array}
\]
mapping into the isotropy group over every point in $S^2 = \mbox{Span}\;\{ e_1 + e_2 + e_3, e_5, e_6\} \cap S^5$.

The first conjugacy class corresponds to a sector diffeomorphic to $Q$.    The next four correspond to sectors diffeomorphic to $S^3$ mod the trivial action of $\Z_3$, already represented by the $\Z$-sectors.  The following three correspond to sectors diffeomorphic to $S^3$ mod the trivial action of $\Z_2$, also diffeomorphic to a $\Z$-sector.  The final three conjugacy classes, however, are comprised of points with isotropy $D_6$ with only the trivial group acting.  These sectors are given by $S^2 = \mbox{Span}\;\{ e_1 + e_2 + e_3, e_5, e_6\} \cap S^5$, and have nonzero Euler-Satake characteristic.  Note that the $\F_2$-sectors correspond to the space of $2$-multisectors (see \cite[page 54]{ademleidaruan}), and that $\F_2$ covers the local groups of $Q$.

We see, then, that although the $\Z$-sectors do not detect any obstruction to the existence of nonvanishing vector fields, the $\F_2$-sectors do.  We note that the obstruction is also detected using $\Gamma = D_6$ and $\Gamma = \Z_2 \oplus \Z_3 \oplus D_6$; the latter of these covers the local groups, while the former does not.

\end{example}


\subsection{Compact Orbifolds With Boundary and Open Suborbifolds of Closed Orbifolds}
\label{subsec-boundarycase}

Although it is likely that the construction of the $\Gamma$-sectors
extends naturally to the case of an orbifold with boundary (whose
orbifold structure is given by a Lie groupoid with $G_0$ and $G_1$
manifolds with boundary), we will not develop the construction in
this case.  Rather, we will use the double orbifold to define them.

Let $Q$ be a compact $n$-dimensional orbifold with boundary (see
\cite{chenruangwt} or \cite{seaton1} for the definition). Form the
double $\widehat{Q}$ (see \cite[Section 3]{paquetteseaton}) and let
$\widehat{\mathcal{G}}$ be an orbifold groupoid for $\widehat{Q}$
with objects $\widehat{G}_0$, arrows $\widehat{G}_1$, source
$\widehat{s}$, target $\widehat{t}$, quotient projection
$\widehat{\sigma}$, etc.  Form the $\Gamma$-sectors
$\widetilde{(\widehat{Q})}_\Gamma$. Treating $Q$ as a subset of
$\widehat{Q}$, we let
\[
    \tilde{Q}_\Gamma = \widetilde{(\widehat{Q})}_\Gamma \cap \widehat{\pi}^{-1}(Q).
\]
Then if $\{ V_x, G_x, \pi_x \}$ is a linear chart at $x \in
\widehat{G}_0$ such that $\widehat{\sigma}(x) = p \in
\partial Q \subseteq Q \subset \widehat{Q}$, a chart for $Q$ can be
taken to be $\{ V_x^+, G_x, \pi_x \}$ where $V_x^+$ is an open
subset of $\R_+^n = \{ (x_1, \ldots ,x_n)  : x_n \geq 0 \}$.
For each $\phi_x \in \widehat{\pi}^{-1}(p)$, a neighborhood of
$\phi_x \in \tilde{Q}_\Gamma$ is covered by the chart with boundary
$\{ V_x^{\langle \phi_x \rangle} \cap V_x^+, C_{G_x}(\phi_x),
\pi_x^{\phi_x} \} = \{(V_x^+)^{\langle\phi_x\rangle},
C_{G_x}(\phi_x), \pi_x^{\phi_x} \}$.  Hence, we see that
$\tilde{Q}_\Gamma$ has the structure of an orbifold with boundary
and $\widetilde{(\widehat{Q})}_\Gamma$ is the double orbifold of
$\tilde{Q}_\Gamma$.

It is easy to see that each $(\phi) \in T_{\widehat{Q}}^\Gamma$ has a representative $\phi_x$ with $\widehat{\pi}(\phi_x) \in Q$, so that we define $T_Q^\Gamma = T_{\widehat{Q}}^\Gamma$.
For each $(\phi) \in T_Q^\Gamma$, let
\[
    \tilde{Q}_{(\phi)} = \widetilde{(\widehat{Q})}_{(\phi)} \cap
    \widehat{\pi}^{-1}(Q).
\]
Then we clearly have
\[
    \tilde{Q}_\Gamma = \bigsqcup\limits_{(\phi) \in T_Q^\Gamma}
    \tilde{Q}_{(\phi)}.
\]
Let $\pi : \tilde{Q}_\Gamma \rightarrow Q$ be defined as the restriction of $\widehat{\pi}$ to $\tilde{Q}_\Gamma \subset \widetilde{(\widehat{Q})}_\Gamma$, and note that the relation $\leq$ defined on $T_{\widehat{Q}}^\Gamma$
coincides with its natural definition on $T_Q^\Gamma$.  In other
words, $(\psi) \leq (\phi)$ as elements of $T_{\widehat{Q}}^\Gamma$
if and only if $\pi\left(\tilde{Q}_{(\psi)}\right) \subseteq
\left(\tilde{Q}_{(\phi)}\right)$.

We also define
\[
    \overline{T_Q^\Gamma} = \left\{ (\phi) \in T_Q^\Gamma :  \partial\left( \tilde{Q}_{(\phi)}\right)  = \emptyset
    \right\}
\]
to be the set of all $\Gamma$ sectors of $Q$ that are closed
orbifolds.  Note that $\Gamma$ covers the local groups of $Q$ if and
only if $\Gamma$ covers the local groups of $\widehat{Q}$.

\begin{theorem}
\label{thrm-boundarycase}

Let $Q$ be a compact orbifold with boundary.  Let $\Gamma$ be a
finitely generated group that covers the local groups of $Q$.  Then
$Q$ admits a smooth nonvanishing vector field if and only if
$\chi_{ES}\left(\tilde{Q}_{(\phi)}\right) = 0$ for each $(\phi) \in
\overline{T_Q^\Gamma}$.

\end{theorem}

Again, we note that no requirement is made of the behavior of the
vector field on the boundary of $Q$.

The proof of this theorem is similar to that of Theorem
\ref{thrm-closedobstruction}.  The primary difference is the
observation that a vector field need not vanish on
$\widehat{\pi}\left(\widetilde{(\widehat{Q})}_{(\phi)}\right)$ for some $(\phi) \in T_Q^\Gamma
\backslash \overline{T}_Q^\Gamma$.  However, since the images of
these sectors intersect the boundary, zeros can be ``pushed off" to
occur outside of $Q$ in $\widehat{Q}$.

\begin{lemma}
\label{lem-boundarynonvanimplieseulerclassvanishes}

Let $Q$ be a compact orbifold with boundary.  If $Q$ admits a
nonvanishing, smooth vector field $X$, then for every finitely
generated $\Gamma$, $\chi_{ES}\left(\tilde{Q}_{(\phi)}\right) = 0$ for
each $(\phi) \in \overline{T}_Q^\Gamma$.

\end{lemma}

\begin{proof}

Suppose $Q$ admits a nonvanishing, smooth vector field $X$.  As $Q$ is closed in $\widehat{Q}$,
we can extend $X$ to a vector field $\widehat{X}$ by Lemma \ref{lem-extendvfields} which, by Lemma \ref{lem-smoothvfldapprox} we can assume is smooth.  By Lemma \ref{lem-bundlesovergammasectors}, $\widehat{X}$ induces a vector field on $\widetilde{(\widehat{Q})}_\Gamma$; let
$\tilde{X}_\Gamma$ denote the restriction of this vector field to $\tilde{Q}_\Gamma$.  Then
$\tilde{X}_\Gamma$ is smooth and nonvanishing.  For each $(\phi) \in
\overline{T}_Q^\Gamma$, we have that $\tilde{Q}_{(\phi)}$ is a
closed orbifold and $\tilde{X}_\Gamma$ is a smooth, nonvanishing
vector field on $\tilde{Q}_{(\phi)}$.  Therefore, the Poincar\'{e}-Hopf Theorem for closed orbifolds in \cite{satake2},
$\chi_{ES}\left(\tilde{Q}_{(\phi)}\right) = 0$.

\end{proof}

Now assume $\Gamma$ covers the local groups of $Q$.

\begin{claim}[Base Case]
\label{clm-boundarybasecase}

Let $Q$ be a compact orbifold with boundary and $\Gamma$ a finitely
generated group that covers the local groups of $Q$. If
$\chi_{ES}\left(\tilde{Q}_{(\phi)}\right) = 0$ for each $(\phi) \in
\overline{T}_Q^\Gamma$, then there is a smooth vector field $X$ on
$Q$ whose restriction to $\pi\left(\tilde{Q}_{(\phi)}\right)$ for
each minimal $(\phi) \in T_Q^\Gamma$ is nonvanishing.

\end{claim}

\begin{proof}

Let $(\phi)$ be a minimal element of $T_Q^\Gamma$.  If $(\phi) \in
\overline{T}_Q^\Gamma$, then as $(\phi)$ is a minimal element of
$T_{\widehat{Q}}^\Gamma$ and $\chi_{ES}\left(\tilde{Q}_{(\phi)}\right)
= 0$, the same technique used in the proof of Claim
\ref{clm-closedbasecase} can be used to construct a nonvanishing,
smooth vector field $X_{(\phi)}$ on $\pi \left( \tilde{Q}_{(\phi)}
\right)$.

If $(\phi) \notin \overline{T}_Q^\Gamma$, then
$\widehat{\pi}\left(\widetilde{\left(\widehat{Q}\right)}_{(\phi)}\right)$ is a
closed manifold by Lemma \ref{lem-minsectorsmanifolds}.  Let
$\widehat{X}_{(\phi)}$ be a smooth vector field on
$\widehat{\pi}\left(\widetilde{\left(\widehat{Q}\right)}_{(\phi)}\right)$ with isolated zeros.  For each
zero of $\widehat{X}_{(\phi)}$ in $p \in \pi\left(\tilde{Q}_{(\phi)}\right)$,
pick a simple smooth curve $c(t)$ in $\widehat{\pi}\left(\widetilde{\left(\widehat{Q}\right)}_{(\phi)}\right)$
with $c(0) = p$ and $c(1) \in
\widehat{Q} \backslash Q$.  Given a tubular neighborhood $W$ of the
image of $c$, the vector field $\widehat{X}_{(\phi)}$ can be smoothly
perturbed on a compact subset of $W$ so that it vanishes only at
$c(1)$.  Applying this to each zero in
$\pi\left(\tilde{Q}_{(\phi)}\right)$ and then letting $X_{(\phi)}$ be the restriction of $\widehat{X}_{(\phi)}$ to $\pi\left(\tilde{Q}_{(\phi)}\right)$, we can assume $X_{(\phi)}$ is
nonvanishing vector field on $\pi\left(\tilde{Q}_{(\phi)}\right)$.

Again, the images of minimal $\Gamma$-sectors are either disjoint or
coincide by Lemma \ref{lem-intersectionssectors}, so we can
construct a nonvanishing vector field on the image of each minimal
$\Gamma$-sector in $Q$. By Lemma \ref{lem-extendvfields}, as the
union of the images of the minimal sectors in $Q$ is closed in $\widehat{Q}$, we can
extend to a vector field $\widehat{X}$ on all of $\widehat{Q}$ which, by Lemma
\ref{lem-smoothvfldapprox} we may assume is smooth.  Then the required $X$ is the restriction of $\widehat{X}$ to $Q$.

\end{proof}

\begin{claim}[Induction Step]
\label{clm-boundaryinductionstep}

Let $Q$ be a compact orbifold with boundary and $\Gamma$ a finitely
generated group that covers the local groups of $Q$. Suppose
$\chi_{ES}\left(\tilde{Q}_{(\phi)}\right) = 0$ for each $(\phi) \in
\overline{T}_Q^\Gamma$.  Let $(\phi) \in T_Q^\Gamma$, and suppose
there is a continuous vector field $X$ on $Q$ that restricts to a
nonvanishing vector field on $\bigsqcup\limits_{(\psi) < (\phi)}
\pi\left(\tilde{Q}_{(\psi)}\right)$.  Then there is a continuous
vector field $Y$ on $Q$ does not vanish on
$\pi\left(\tilde{Q}_{(\phi)}\right)$ and coincides with $X$ on each
$\pi\left(\tilde{Q}_{(\psi)}\right)$ with $(\psi) < (\phi)$.

\end{claim}

\begin{proof}

If $(\phi) \in \overline{T}_Q^\Gamma$ then
$\pi\left(\tilde{Q}_{(\phi)}\right)$ is closed, so the proof is
identical to that of Claim \ref{clm-closedinductionstep}.  On the
other hand, if $(\phi) \notin \overline{T}_Q^\Gamma$, then the
set $\pi \left(\widetilde{\left(\widehat{Q}\right)}_{(\phi)}\right)
\backslash \bigcup\limits_{(\psi) < (\phi)}
\pi\left(\widetilde{\left(\widehat{Q}\right)}_{(\psi)}\right)$ is a manifold.  Therefore, $X|_{\pi
\left(\widetilde{\left(\widehat{Q}\right)}_{(\phi)}\right)}$ can be
continuously perturbed away from each of the
$\pi\left(\widetilde{\left(\widehat{Q}\right)}_{(\psi)}\right)$ with
$(\psi) < (\phi)$ so that it vanishes only on $\widehat{Q}
\backslash Q$.  The resulting vector field can be extended by Lemma
\ref{lem-extendvfields} to a continuous vector field $\widehat{Y}$ on $\widehat{Q}$; the required vector field $Y$ is the restriction $\widehat{Y}$ to $Q$.

\end{proof}

\begin{proof}[Proof of Theorem \ref{thrm-boundarycase}]

Given Claims \ref{clm-boundarybasecase} and Claim
\ref{clm-boundaryinductionstep}, the proof of Theorem
\ref{thrm-boundarycase} is identical to the proof of Theorem
\ref{thrm-closedobstruction}.

\end{proof}

Techniques almost identical to those above can be used to prove the following.

\begin{theorem}
\label{thrm-opencase}

Let $Q$ be an open suborbifold of the closed orbifold $R$. Let
$\Gamma$ be a finitely generated group that covers the local groups
of $Q$. Then $Q$ admits a smooth nonvanishing vector field if and
only if $\chi_{ES}\left(\tilde{Q}_{(\phi)}\right) = 0$ for each $(\phi)
\in T_Q^\Gamma$ such that $\tilde{Q}_{(\phi_x)}$ is a closed orbifold.

\end{theorem}

For this case, we note that each $(\phi) \in T_Q^\Gamma$ determines
an $\approx$-class in $T_R^\Gamma$.  The correspondence
$(\phi)_\approx^Q \mapsto (\phi)_\approx^R$ is neither surjective
nor injective.  Rather than use this correspondence, we apply the induction in the proof of
Theorem \ref{thrm-closedobstruction} to the sectors of
$T_R^\Gamma$ using the techniques in the proof of Theorem
\ref{thrm-boundarycase} when $\pi\left(\tilde{Q}_{(\phi)}\right)$ is
not completely contained in $R$; i.e. when $\tilde{Q}_{(\phi)}$ is
not closed.  With this observation, the argument extends directly.


\bibliographystyle{amsplain}

\end{document}